\documentclass[12pt, leqno]{article}
\usepackage{amsfonts}
\usepackage{amsmath}
\usepackage{amssymb}
\usepackage{amscd}
\begin{document}
\renewcommand{\thefootnote}{\fnsymbol{footnote}}
\begin{center}
{\bf \Large $p$-Johnson homomorphisms and pro-$p$ groups}\footnote[0]{
$\;\;\;\;$ Key words: arithmetic topology, pro-$p$ groups, $p$-Johnson homomorphisms,  Zassenhaus filtration, non-Abelian Iwasawa theory, Massey products.
}
\end{center} 

\vspace{0.08cm}

\begin{center}
Masanori Morishita  and Yuji Terashima 
\end{center}

\begin{center}
{\it Dedicated to Professor Kazuya Kato}
\end{center}

\vspace{0.05cm}

{\small {\bf Abstract.} We propose an approach  to study non-Abelian Iwasawa theory, using the idea of Johnson homomorphisms in low dimensional topology. We introduce arithmetic analogues of  Johnson homomorphisms/maps, called the $p$-Johnson homomorphisms/maps, associated to the Zassenhaus filtration of a pro-$p$ Galois group over a $\mathbb{Z}_p$-extension of a number field. We give their cohomological interpretation in terms of Massey products in Galois cohomology.}

\vspace{0.08cm}

\begin{center}
{\bf 1. Introduction}
\end{center}

Let $p$ be an odd prime number, and let $\mu_{p^n}$ denote the group of $p^n$-th roots of unity for a positive integer $n$ and we set $\displaystyle{\mu_{p^{\infty}} := \cup_{n\geq 1} \mu_{p^n}}$.  We let $k_{\infty} := \mathbb{Q}(\mu_{p^{\infty}})$  and $\tilde{k}$ the maximal pro-$p$ extension of $k_{\infty}$ which is unramified outside $p$. We let $\Gamma_p := {\rm Gal}(k_\infty/\mathbb{Q})$ and $F_p := {\rm Gal}(\tilde{k}/k_\infty)$,
the Galois groups of the extensions $k_\infty/\mathbb{Q}$ and $\tilde{k}/k_\infty$, respectively. Classical Iwasawa theory then deals with the action of $\Gamma_p$ on the Abelianization $H_1(F_p,\mathbb{Z}_p)$ of $F_p$ ([Iw1]).
A basic problem of non-Abelian Iwasawa theory, with which we are concerned in this paper,  is to study the conjugate action of $\Gamma_p$ on $F_p$ itself. In terms of schemes, one has the tower of \'{e}tale pro-finite covers
 
{\small $$ \tilde{X}_p := {\rm Spec}({\cal O}_{\tilde{k}}[1/p]) \rightarrow X_p^{\infty} := {\rm Spec}({\cal O}_{k_\infty}[1/p]) \rightarrow X_p := {\rm Spec}(\mathbb{Z}[1/p]), \leqno{(1.1)}$$}

where ${\cal O}_{k_\infty}$ and ${\cal O}_{\tilde{k}}$ denote the rings of integers of $k_{\infty}$ and $\tilde{k}$, respectively, and the Galois groups
$$ \Gamma_p = {\rm Gal}(X_p^{\infty}/X_p), \;\; F_p = {\rm Gal}(\tilde{X}_p/X_p^{\infty}) = \pi_1^{{\footnotesize {\rm pro}-p}}(X_p^{\infty}),  \leqno{(1.2)}$$
where $\pi_1^{{\footnotesize {\rm pro-}p}}$ stands for the maximal pro-$p$ quotient of the \'{e}tale fundamental group. So the problem is to study the monodromy action of $\Gamma_p$ on the arithmetic pro-$p$ fundamental group $F_p$.
\\

Now let us recall  the analogy between a prime and a knot \\
$$
\mbox{
\begin{tabular}{|c|c|}
\hline
prime& knot   \\
${\rm Spec}(\mathbb{F}_p) = K(\hat{\mathbb{Z}},1) \hookrightarrow \overline{{\rm Spec}(\mathbb{Z})}$ & $ S^1 = K(\mathbb{Z},1) \hookrightarrow S^3$  \\
\hline
\end{tabular}
}
\leqno{(1.3)}
$$
\\
Here  $K(*,1)$ stands for  the Eilenberg-MacLane space and $\overline{{\rm Spec}(\mathbb{Z})} := {\rm Spec}(\mathbb{Z}) \cup \{ \infty\}$, $\infty$ being the infinite prime of $\mathbb{Q}$ which may be seen as an analogue of the end of $\mathbb{R}^3$ ([De]). This analogy (1.3) opens a research area, called {\it arithmetic topology}, which studies systematically further analogies between number theory and 3-dimentional topology  ([Ms2]). In particular, there are known intimate analogies between Iwasawa theory and Alexander-Fox theory ([Ma], [Ms2; Chap. $9 \sim 12$]).\\
 
Arithmetic topology suggests that topological counterparts of (1.1) and (1.2) may be the tower of covers
$$ \tilde{X}_{\cal K} \rightarrow X_{\cal K}^{\infty} \rightarrow X_{\cal K} := S^3 \setminus {\cal K}, $$
for a knot ${\cal K}$ in $S^3$, where $X_{\cal K}^\infty$ and $\tilde{X}_{\cal K}$ denote the infinite cyclic cover and the universal cover of the knot complement $X_{\cal K}$, respectively, and the Galois groups
$$ \Gamma_{\cal K} := {\rm Gal}(X_{\cal K}^\infty/X_{\cal K}),\;\; F_{\cal K} := {\rm Gal}(\tilde{X}_{\cal K}/X_{\cal K}^\infty) = \pi_1(X_{\cal K}^\infty), $$
and we have the conjugate action of $\Gamma_{\cal K}$ on $F_{\cal K}$.

To push our idea further,  suppose that ${\cal K}$ is a fibered knot so that $X_{\cal K}$ is a mapping torus 
of the monodromy $\phi : S \rightarrow S$, $S$ being the Seifert surface of ${\cal K}$. Then $F_{\cal K} = \pi_1(S)$ and the conjugate action of $\Gamma_{\cal K}$ on $F_{\cal K}$ is nothing but the monodromy action induced by $\phi$ on $F_{\cal K}$
$$ \phi_* \; : \; \Gamma_{\cal K} \longrightarrow {\rm Aut}(F_{\cal K}). \leqno{(1.4)}$$
Note here that the monodromy $\phi$ may be regarded as a mapping class of the surface $S$. Thus the action (1.4) can be studied by means of  the Johnson homomorphisms/maps, associated to the lower central series of $F_{\cal K}$,  defined on a certain filtration of the mapping class group for the surface $S$  ([J], [Ki], [Mt]) or, more generally, on the automorphism group ${\rm Aut}(F_{\cal K})$  ([Kw], [Sa]). \\

In this paper, we regard the action of $\Gamma_p$ on $F_p$ as an arithmetic analogue of the monodromy action (1.4) and propose an approach to study non-Abelian Iwasawa theory by introducing arithmetic analogues of the Johnson homomorphisms/maps,  called the $p$-Johnson homomorphisms/maps, associated to the Zassenhaus filtration of $F_p$, defined on a certain filtration of the automorphism group ${\rm Aut}(F_p)$.   For this, we lay a foundation of a general theory of $p$-Johnson homomorphisms/maps in the context of pro-$p$ group. \\

We note that our viewpoint and approach differs from what is called ``non-commutative Iwasawa theory" (cf. [CFKSV], [Kt; 3]). The works by  M. Ozaki ([O]) and R. Sharifi ([Sh]) are related to ours  (see Remark 3.2.7), however, our approach is different from theirs and closer to geometric topology.\\

Here is the content of this paper. In Section 2,  we give a general theory of $p$-Johnson homomorphisms in the context of pro-$p$ groups.  We use the Zassenhaus filtration of a finitely generated pro-$p$ group $G$ in order to  introduce the  $p$-Johnson homomorphisms, defined on a certain filtration of the automorphism group of $G$. In Section 3, we give a framework to study  non-Abelian Iwasawa theory by means of the $p$-Johnson homomorphisms.  In Section 4, we give a theory of Johnson maps for a free pro-$p$ group $F$ by  extending  the $p$-Johnson homomorphisms in Section 2 to maps, called the $p$-Johnson maps,  defined on the automorphism group ${\rm Aut}(F)$ itself. In Section 5, we give a cohomological interpretation of the $p$-Johnson homomorphisms in terms of Massey products in Galois cohomology.\\
\\
{\it Notation.} For subgroup $A, B$ of a group $G$, $[A,B]$ stands for the subgroup of $G$ generated by $[a,b] := aba^{-1}b^{-1}$ for all $a \in A, b \in B$.\\

\begin{center}
{\bf 2.  Zassenhaus filtration and  $p$-Johnson homomorphisms for a pro-$p$ group.} 
\end{center} 

In this section, we give a general theory of $p$-Johnson homomorphisms for pro-$p$ groups. We associate to the Zassenhaus filtration of a finitely generated pro-$p$ group $G$ a certain filtration on the automorphism group ${\rm Aut}(G)$ of $G$, and introduce the $p$-Johnson homomorphisms defined on each term of the filtration of ${\rm Aut}(G)$. 

Throughout this section, let $p$ be a fixed prime number and $G$ a finitely generated pro-$p$ group. For general properties of pro-$p$ groups, we consult [Ko] and  [DDMS].
\\
\\
{\bf 2.1. Zassenhaus filtration and the associated Lie algebra.}  Let $\mathbb{F}_p[[G]]$ be the complete group algebra of $G$ over $\mathbb{F}_p = \mathbb{Z}/p\mathbb{Z}$ with the augmentation ideal $I_G := {\rm Ker}(\epsilon_{\mathbb{F}_p[[G]]})$, where  $\epsilon_{\mathbb{F}_p[[G]]} : \mathbb{F}_p[[G]] \rightarrow \mathbb{F}_p$ is the augmentation homomorphism ([Ko; 7.1]). For each positive integer $n$, we define the normal subgroup $G_n$ of $G$ by
$$ G_n := \{ g \in G \; | \; g-1 \in I_G^n \}. \leqno{(2.1.1)}$$
The descending series $\{G_n\}_{n\geq 1}$ is called the {\it Zassenhaus filtration} of $G$ ([Ko; 7.4]). The family $\{ G_n\}_{n \geq 1}$ forms a full system of neighborhoods of the identity $1$ in $G$ and satisfies the following properties
$$ \begin{array}{ll}
(2.1.2)& \;\;\;\;\; \;\;\;\;\;\;\;\;\;\; (G_i)^p \subset G_{pi} \;\;\;\;\;\; (i \geq 1).\\
(2.1.3)& \;\;\;\;\; \;\;\;\;\;\;\;\;\;\; [G_i, G_j] \subset G_{i+j} \;\;\;\;\; \; (i, j \geq 1).\\
\end{array}
$$
We recall the fact that the abstract commutator subgroup of a finitely generated pro-$p$ group is closed ([DDMS;1.19]). 

The Zassenhaus filtration is in fact the fastest descending series of $G$ having the properties (2.1.2) and (2.1.3). Namely, it is shown by Jennings' theorem and an inverse limit argument that  we have the following inductive description of $G_n$:
$$ G_n = (G_{[n/p]})^p \prod_{i+j = n} [G_i,G_j] \;\;\;\;\; (n\geq 2), \leqno{(2.1.4)}$$
where $[n/p]$ stands for the least integer $m$ such that $mp \geq n$. ([DDMS; 12.9]). 

We note by (2.1.3) that elements of $G_i/G_{i+j}$ and $G_j/G_{i+j}$ commute, in particular, $G_n/G_{n+1}$ is central in $G/G_{n+1}$. The 2nd term $G_2$ is the Frattini subgroup $G^p[G,G]$ of $G$ and  we denote by $H$ the Frattini quotient 
$$ H := G/G_2 = G/G^p[G,G] = H_1(G,\mathbb{F}_p). \leqno{(2.1.5)}$$
For $g \in G$, we write $[g]$ for the image of $g$ in $H$: $[g] := g$ mod $G_2$. We note that each $G_n$ is a finitely generated pro-$p$ group ([DDMS; 1.7, 1.14]). \\

For each $n \geq 1$, we let
$$ {\rm gr}_n(G) := G_n/G_{n+1}, $$ 
which is a finite dimensional $\mathbb{F}_p$-vector space.  The graded $\mathbb{F}_p$-vector space 
$$ \displaystyle{{\rm gr}(G) := \bigoplus_{n\geq 1} {\rm gr}_n(G)} \leqno{(2.1.6)}$$
has a natural structure of a graded Lie algebra over $\mathbb{F}_p$  by (2.1.3). Here, for $a = g \; {\rm mod}\; G_{i+1}, b = h \; {\rm mod}\; G_{j+1} \;\; (g \in G_i, h \in G_j)$, the Lie bracket is defined by
$$ [a, b]_{{\rm gr}(G)} := [g,h] \; {\rm mod}\; G_{i+j+1}. $$
Further, by (2.1.2) again, ${\rm gr}(G)$ has the operation $[p]$ defined by, for $a = g \; {\rm mod}\; G_{n+1} \in {\rm gr}_n(G)$, 
$$ [p](a) := g^p \; {\rm mod}\; G_{pn+1},$$
which makes ${\rm gr}(G)$ a restricted Lie algebra over $\mathbb{F}_p$ ([DDMS; 12.1]). 

The {\it restricted universal enveloping algebra} (abbreviated to {\it universal envelope}) $U({\rm gr}(G))$ of  ${\rm gr}(G)$ is given as follows. For each $m \geq 0$, we let
$$ {\rm gr}_m(\mathbb{F}_p[[G]]) :=  I_G^m/I_G^{m+1}.$$
and consider the graded associative algebra over $\mathbb{F}_p$:
$$ \displaystyle{ {\rm gr}(\mathbb{F}_p[[G]]) :=  \bigoplus_{m\geq 0} {\rm gr}_m(\mathbb{F}_p[[G]]).}$$
For each $m\geq 1$, we have an injective $\mathbb{F}_p$-linear map
$$ \theta_m \; :\: {\rm gr}_m(G) \longrightarrow  {\rm gr}_m(\mathbb{F}_p[[G]])$$
defined by
$$ \theta_m( g \; {\rm mod}\; G_{m+1}) := g-1 \; {\rm mod}\; I_G^{m+1}\;\; \mbox{for} \; g \in G_m. $$
 Putting all $\theta_m$ together over $m \geq 1$, we have an injective graded Lie algebra homomorphism over $\mathbb{F}_p$
$$ {\rm gr}(\theta) := \bigoplus_{m \geq 1} \theta_m \; : \; {\rm gr}(G) \longrightarrow {\rm gr}(\mathbb{F}_p[[G]]).$$
Then $({\rm gr}(\mathbb{F}_p[[G]]), {\rm gr}(\theta))$ is the universal envelope of ${\rm gr}(G)$ ([DDMS; 12.8]):
$$ U{\rm gr}(G)  = {\rm gr}(\mathbb{F}_p[[G]]). \leqno{(2.1.7)}$$
\\
{\bf 2.2. The automorphism group and $p$-Johnson homomorphisms.} Let ${\rm Aut}(G)$ denote the group of continuous automorphisms of a finitely generated pro-$p$ group $G$. We note that any abstract group homomorphism between finitely generated pro-$p$ groups
is always continuous  and so ${\rm Aut}(G)$ is same as the group of automorphisms  of $G$ (as an abstract group) ([DDMS; 1.21]). We also note that every term $G_n$ of the Zassenhaus filtration of $G$ is a characteristic subgroup of $G$, namely, invariant under the action of ${\rm Aut}(G)$. 

Since any automorphism $\phi$ of $G$ induces an automorphism $[\phi]_m$ of $G/G_{m+1}$ for each integer $m \geq 0$, we have the group homomorphism 
$$ [\;\;]_m \; : \; {\rm Aut}(G) \longrightarrow {\rm Aut}(G/G_{m+1}). \leqno{(2.2.1)}$$
We then define the normal subgroup ${\rm A}_G(m)$ of ${\rm Aut}(G)$ by
$$ \begin{array}{ll} {\rm A}_G(m) & := {\rm Ker}([\;\;]_m)  \\
                    & = \{ \phi \in {\rm Aut}(G) \, | \, \phi(g)g^{-1} \in G_{m+1} \} \;\;\;\; (m \geq 0). 
\end{array}
 \leqno{(2.2.2)} 
$$
We call the resulting descending series $\{ {\rm A}_G(m) \}_{m \geq 0}$ the {\it Andreadakis-Johnson filtration} of ${\rm Aut}(G)$ associated to the Zassenhaus filtration of $G$  (cf [A], [Sa]). In particular, we set simply  $[\phi] := [\phi]_1$ for $\phi \in {\rm Aut}(G)$ and the 1st term ${\rm A}_G(1)$  is called the {\it induced automorphism group} of $G$ and denoted by
${\rm IA}(G)$:
$$ {\rm IA}(G) := {\rm Ker}([\;\;]: {\rm Aut}(G) \longrightarrow {\rm GL}(H)), \leqno{(2.2.3)}$$
where ${\rm GL}(H)$ denotes the group of $\mathbb{F}_p$-linear automorphisms of $H = G/G_2$. 

The family $\{ {\rm A}_G(m) \}_{m\geq 0}$ forms a full system of neighborhood  of the identity ${\rm id}_G$ in ${\rm Aut}(G)$ and it can be shown that ${\rm Aut}(G)$ is a pro-finite group
 and ${\rm IA}(G)$ is a pro-$p$ group ([DDMS; 5.3, 5.5]). So ${\rm Aut}(G)$ is virtually a pro-$p$ group.\\
 
The next Lemma will play a basic role to introduce the $p$-Johnson homomorphisms. \\
\\
{\bf Lemma 2.2.4.} {\it For $\phi \in {\rm A}_G(m)$ $(m\geq 0)$ and $g \in G_n$ $(n \geq 1)$, we have}
$$ \phi(g)g^{-1} \in  G_{m+n}.$$
{\it Proof.}  We fix $m$ and prove the assertion by induction on $n$. For $n =1$, the assertion $\phi(g)g^{-1} \in G_{m+1}$ is true by definition (2.2.2) of ${\rm A}_G(m)$. Assume that
$$ \phi(g)g^{-1} \in G_{m+i} \; \mbox{if} \; g \in G_i \; \mbox{and} \; 1\leq i \leq n. \leqno{(2.2.4.1)}$$ 
By (2.1.4), we have
$$ \displaystyle{G_{n+1} = (G_{[(n+1)/p]})^p \prod_{i+j=n+1} [G_i,G_j].}$$
Since $G_{n+1}/(\prod_{i+j=n+1}[G_i,G_j])$ is Abelian, we have
$$\displaystyle{ G_{n+1} = \{ a^p \, | \, a \in G_{[(n+1)/p]} \} \prod_{i+j=n+1}[G_i,G_j]}$$
and so any element $g$ of $G_{n+1}$ can be written in the form
$$ g = a^p [b_1,c_1]^{e_1}\cdots [b_q,c_q]^{e_q},$$
where $a \in G_{[(n+1)/p]}$ and for each $s$ ($1\leq s \leq q$) there are $i, j$ ($i+j=n+1$) such that $b_s \in G_i, c_s \in G_j$. Since we have
$$ \phi(g)g^{-1} = \phi(a)^p\phi([b_1,c_1])^{e_1}\cdots \phi([b_q,c_q])^{e_q}[b_q,c_q]^{-e_q}\cdots [b_1,c_1]^{-e_1}a^{-p},$$
it suffices to show that
$$ \left\{
\begin{array}{l}
(2.2.4.2) \;\;\phi([b,c])[b,c]^{-1} \in G_{m+n+1} \; \mbox{if} \; b \in G_i, c \in G_j \; \mbox{and}\; i+j=n+1,\\
(2.2.4.3)\;\; \phi(a)^p a^{-p} \in G_{m+n+1} \; \mbox{if} \; a \in G_{[(n+1)/p]}.
\end{array}
\right.
$$
\\
(2.2.4.2). For simplicity, we shall use the notation: $[\psi, x] := \psi(x)x^{-1}$ and $[x,\psi] := x \psi(x)^{-1}$ for $x \in G$ and $\psi \in {\rm Aut}(G)$. 
By  the ``three subgroup lemma" and the induction hypothesis (2.2.4.1), we have
$$ \begin{array}{ll}
\phi([b,c])[b,c]^{-1} & = [\phi, [b,c]] \\
                                  & \in [\phi, [G_i,G_j]]\\
                                  & \subset [[\phi,G_i],G_j][[G_j,\phi],G_i]\\
                                  & \subset [G_{m+i},G_j][G_{m+j},G_i]\\
                                  & = G_{m+i+j} = G_{m+n+1}.
\end{array}
$$
(2.2.4.3). Let $t := [(n+1)/p]$ so that $pt \geq n+1$. By (2.1.1) and the induction hypothesis (2.2.4.1), we have
$$ \phi(a) - a = (\phi(a)a^{-1}-1)a \in I_G^{t+m}.$$
Therefore we have
$$ \begin{array}{ll}
\phi(a)^p a^{-p} -1  &= (\phi(a)^p -a^p)a^{-p}\\
                                & = (\phi(a)-a)^p a^{-p}\\
                                & \in I^{p(t+m)} \subset I^{m+n+1}.
\end{array}
$$
 Hence $\phi(a)^p a^{-p} \in G_{m+n+1}$ by (2.1.1). $\;\; \Box$\\
\\                               
 Lemma 2.2.4 yields the following properties of the Andreadakis-Johnson filtration $\{ {\rm A}_G(m)\}_{m\geq 0}$.\\
 \\
 {\bf Proposition 2.2.5.}  {\it We have}
 $$ \begin{array}{ll} (1) & \;\;\;\;\; \;\;\;\;\;\;\;\;   [{\rm A}_G(i), {\rm A}_G(j)] \subset {\rm A}_G(i+j) \;\; \mbox{for}\, i, j \geq 0.\\
 (2)  & \;\;\;\;\; \;\;\;\;\;\;\;\;\;\;  {\rm A}_G(m)^p \subset {\rm A}_G(m+1) \;\; \mbox{if}\, m \geq 1.
 \end{array}
 $$
 {\it Proof.} (1) We use the same notation as in the proof of (2.2.4.2). By Lemma 2.2.4, we have
 $$ \begin{array}{l}
 [[{\rm A}_G(j),G],{\rm A}_G(i)] \subset [G_{j+1}, {\rm A}_G(i)] \subset G_{i+j+1},\\
 
 [[G,{\rm A}_G(i)], {\rm A}_G(j)] \subset [G_{i+1}, {\rm A}_G(j)] \subset G_{i+j+1}.
 \end{array}
 $$
 By the three subgroup lemma, we have
 $$ [[{\rm A}_G(i), {\rm A}_G(j)], G] \subset [{\rm A}_G(j),G],{\rm A}_G(i)] [[G,{\rm A}_G(i)], {\rm A}_G(j)] \subset G_{i+j+1}.$$
 By definition (2.2.2), we obtain
 $$ [{\rm A}_G(i), {\rm A}_G(j)] \subset {\rm A}_G(i+j).$$
 (2) Let $ g \in G$ and $\phi \in {\rm A}_G(m)$. We shall show that for any integer  $d \geq 1$, 
 $$ \phi^d(g)g^{-1} \equiv (\phi(g)g^{-1})^d \; \mbox{mod}\; G_{2m+1}, \leqno{(2.2.5.1)}$$
 from which the assertion follows. In fact,  let $d = p$ in (2.2.5.1). Then  $(\phi(g)g^{-1})^p \in G_{p(m+1)}$ by (2.1.2), and $G_{2m+1} \subset
 G_{m+2}$ because $m \geq 1$. So $\phi^p(g)g^{-1} \in G_{m+2}$ and hence $\phi^p \in {\rm A}_G(m+1)$.
 
 We prove (2.2.5.1) by induction on $d$. For $d = 1$, it is obviously true. Suppose $\phi^d(g)g^{-1} \equiv (\phi(g)g^{-1})^d \; \mbox{mod}\; G_{2m+1}$.
 Note that $\phi^d(g)g^{-1} \in G_{m+1},$ since $(\phi(g)g^{-1})^d \in G_{m+1}$. Then we have 
$$ \begin{array}{ll}
\phi^{d+1}(g) g^{-1} (\phi(g)g^{-1})^{-(d+1)} & = \phi^{d+1}(g) \phi(g)^{-1} \phi(g) g^{-1} (\phi(g)g^{-1})^{-(d+1)}\\
                                                                             & = \phi(\phi^d(g) g^{-1}) (\phi(g)g^{-1})^{-d} \\
                                                                             & \equiv \phi(\phi^d(g) g^{-1}) (\phi^d(g)g^{-1})^{-1} \; {\rm mod} \; G_{2m+1}. 
\end{array}
$$
Since $ \phi(\phi^d(g) g^{-1}) (\phi^d(g)g^{-1})^{-1} \in G_{2m+1}$ by Lemma 2.2.4,  $\phi^{d+1}(g)g^{-1}   \equiv (\phi(g)g^{-1})^{d+1}$ mod $G_{2m+1}$ and hence the induction holds. $\;\;\Box$  \\

Now we are going to introduce the $p$-Johnson homomorphisms.  Let $\phi \in {\rm A}_G(m)$ ($m \geq 0$). For $g \in G$, we have $\phi(g)g^{-1} \in G_{m+1}$. Then we see that $\phi(g)g^{-1}$ mod $G_{m+2} \in {\rm gr}_{m+1}(G)$ depends only on the class $[g] \in H$. In fact, for $g' = g g_2$ with $g_2 \in G_2$, we have
$$ \phi(g')g'^{-1} = \phi(g) \phi(g_2) g_2^{-1} g^{-1} \equiv \phi(g)g^{-1}  \; {\rm mod}\; G_{m+2},$$
since $\phi(g_2)g_2^{-1} \in G_{m+2}$ by Lemma 2.2.4. Thus we have a map
$$ \tau_m(\phi) \; : \; H \longrightarrow {\rm gr}_{m+1}(G) $$
defined by
$$ \tau_m(\phi)(h) := \phi(g)g^{-1} \; {\rm mod} \; G_{m+2} \;\; (h = [g]). \leqno{(2.2.6)}$$
\\
{\bf Lemma 2.2.7.} {\it For $\phi \in {\rm A}_G(m)$ $(m \geq 0)$, the map $\tau_m(\phi)$ is $\mathbb{F}_p$-linear.}\\
\\
{\it Proof.} Let $h = [g], h' = [g']$ and $c \in \mathbb{F}_p$. Using the property that $G_{m+1}/G_{m+2}$ is central in $G/G_{m+2}$, we have
$$ \begin{array}{ll}
\tau_m(\phi)(h+h') & = \tau_m(\phi)([gg'])\\
                             & = \phi(g g')(g g')^{-1} \; {\rm mod} \; G_{m+2}\\
                            & =  \phi(g) \phi(g')g'^{-1}  g^{-1} \; {\rm mod} \; G_{m+2}\\
                             & = (\phi(g)g^{-1})(\phi(g')g'^{-1}) \; {\rm mod} \; G_{m+2}\\
                             & = \tau_m(\phi)(h) + \tau_m(\phi)(h'),
 \end{array}
$$                            
and
$$ \begin{array}{ll}
\tau_m(\phi)(ch) & = \tau_m(\phi)([g^c]) \\
                          & = \phi(g^c)g^{-c} \; {\rm mod}\; G_{m+2}\\
                          & = (\phi(g)g^{-1})^c \; {\rm mod}\; G_{m+2}\\
                          &  = c \tau_m(\phi)(h).  \;\;\;\;\; \Box
\end{array}
$$                          
\\
Let ${\rm Hom}_{\mathbb{F}_p}(H,{\rm gr}_{m+1}(G))$ denote the group of $\mathbb{F}_p$-linear maps $H \rightarrow {\rm gr}_{m+1}(G)$. By Lemma 2.2.7, we have the map
$$ \tau_m \; :\; {\rm A}_G(m) \longrightarrow {\rm Hom}_{\mathbb{F}_p}(H,{\rm gr}_{m+1}(G)).  $$
For $m=0$, we easily see by (2.2.6) that $\tau_0(\phi) = [\phi] - {\rm id}_H$ for $\phi \in {\rm Aut}(G)$.\\
\\
{\bf Theorem 2.2.8.} {\it  For $m \geq 1$, the map $\tau_m$ is a group homomorphism and  its kernel is ${\rm A}_G(m+1)$.}\\
\\
{\it Proof.} Let $\phi_1, \phi_2 \in {\rm A}_G(m)$.  For any $g \in G$, we have
$$ \begin{array}{ll}
\tau_m(\phi_1 \phi_2)([g]) & = \phi_1(\phi_2(g))g^{-1} \; {\rm mod}\; G_{m+2}\\
                                       & = \phi_1(\phi_2(g) g^{-1})\cdot \phi_1(g)g^{-1} \; {\rm mod}\; G_{m+2}.
\end{array}
$$
Since $\phi_2(g)g^{-1} \in G_{m+1}$, $\phi_1(\phi_2(g)g^{-1}) \equiv \phi_2(g)g^{-1}$ mod $G_{2m+1}$ by Lemma 2.2.4.  Since $G_{2m+1} \subset G_{m+2}$ by $m\geq 1$, we have
$$ \begin{array}{ll}
\tau_m(\phi_1 \phi_2)([g]) & = \phi_1(g) g^{-1}\cdot \phi_2(g)g^{-1} \; {\rm mod}\; G_{m+2}\\
                                          & = (\tau_m(\phi_1) + \tau_m(\phi_2))([g])
                                       \end{array}
                                       $$
for any $g \in G$. Hence the former assertion is proved. The latter assertion on ${\rm Ker}(\tau_m)$ is obvious by definition (2.2.6). $\;\; \Box$
\\
\\
The homomorphism $\tau_m :  {\rm A}_G(m) \rightarrow {\rm Hom}_{\mathbb{F}_p}(H,{\rm gr}_{m+1}G))$ ($m \geq 1$) or the induced injective homomorphism
$$ \overline{\tau}_m \, : \, {\rm gr}_m({\rm A}_G) := {\rm A}_G(m)/{\rm A}_G(m+1)  \hookrightarrow {\rm Hom}_{\mathbb{F}_p}(H,{\rm gr}_{m+1}(G)) \;\;\;\; (m\geq 1)$$
is called the $m$-th {\it $p$-Johnson homomorphism}.\\

We give some properties of the $p$-Johnson homomorphisms. Firstly, we note that the group ${\rm Aut}(G)$ acts on both ${\rm A}_G(m)$ and ${\rm Hom}_{\mathbb{F}_p}(H, {\rm gr}_{m+1}(G))$ by the following rules, respectively:
$$\left\{ 
\begin{array}{l}
\psi.\phi :=   \psi \circ \phi \circ \psi^{-1} \;\; (\psi \in {\rm Aut}(G), \phi \in {\rm A}_G(m)),\\
(\psi.\eta)(h) := \psi(\eta([\psi]^{-1}(h))) \;\; (\psi \in {\rm Aut}(G), \eta \in {\rm Hom}_{\mathbb{F}_p}(H, {\rm gr}_{m+1}(G)), h \in H).
\end{array}
\right.
$$
Then we have the following \\
\\
{\bf Proposition 2.2.9.} {\it The $p$-Johnson homomorphism $\tau_m$ $(\mbox{resp.}\;\overline{\tau}_m)$ is ${\rm Aut}(G)$-equivariant $(\mbox{resp.} \;{\rm Aut}(G)/{\rm IA}(G)$-equivariant$)$.}\\
\\
{\it Proof.} Let $\psi \in {\rm Aut}(G)$ and $\phi \in {\rm A}_G(m)$. Then we have, for any $g \in G$,
$$ \begin{array}{ll}
\tau_m(\psi.\phi)([g]) & = \tau_m(\psi \circ \phi \circ \psi^{-1})([g])\\
                                 & = (\psi \circ \phi \circ \psi^{-1})(g) g^{-1} \; \mbox{mod}\; G_{m+2}.
\end{array}
$$
On the other hand, we have, for any $g \in G$, 
$$ \begin{array}{ll}
(\psi.\tau_m(\phi))([g]) & = \psi(\tau_m(\phi)([\psi]^{-1}([g])))\\
                                 & = \psi(\tau_m(\phi))([\psi^{-1}(g)]))\\
                                 & = \psi(\phi(\psi^{-1}(g)) (\psi^{-1}(g))^{-1}) \; {\rm mod} \; G_{m+2}\\
                                 & =  (\psi \circ \phi \circ \psi^{-1})(g) g^{-1} \; \mbox{mod}\; G_{m+2}.
\end{array}
$$         
Hence $\tau_m$ is ${\rm Aut}(G)$-equivariant. As for $\overline{\tau}_m$, it suffices to note that ${\rm IA}(G)$ acts trivially on ${\rm gr}_m({\rm A}_G) = {\rm A}_G(m)/{\rm A}_G(m+1)$ 
by Proposition 2.2.5 (1) and on ${\rm Hom}_{\mathbb{F}_p}(H, {\rm gr}_{m+1}(G))$ by (2.2.3) and Lemma 2.2.4. $\;\; \Box$\\
\\
Next we compute the $p$-Johnson homomorphism on inner automorphisms. Let $ {\rm Inn} : G \rightarrow {\rm Aut}(G)$ be the homomorphism defined by 
$$ {\rm Inn}(x)(g) := x g x^{-1} \;\; (x, g \in G).$$
The image ${\rm Inn}(G)$ is a normal subgroup of ${\rm Aut}(G)$ and called the group of {\it inner automorphisms} of $G$.
\\
\\
{\bf Proposition 2.2.10.} {\it Let $m \geq 1$ and $x \in G_{m}$. Then  we have}
$$ {\rm Inn}(x) \in {\rm A}_G(m)$$
and  
$$ \tau_m({\rm Inn}(x))([g]) = [x,g] \; {\rm mod}\; G_{m+2} \;\; \;\; (g \in G).$$
{\it Proof.} For $x \in G_{m}$ and $g \in G$, we have 
$${\rm Inn}(x)(g)g^{-1} = [x,g] \in G_{m+1},$$
from which the assertions follow. $\;\; \Box$\\
\\
Finally we compute the $p$-Johnson homomorphisms on commutators of automorphisms.\\
\\
{\bf Lemma 2.2.11.} {\it For $\psi \in {\rm A}_G(i), \phi \in {\rm A}_G(j)$ $(k, m\geq 0)$ and $g \in G$, we have, in ${\rm gr}_{i+j+1}(G),$}
$$ \begin{array}{l} \tau_{i+j}([\psi, \phi])([g]) \\
\;\;\;\;\; = \psi(\phi(g)g^{-1})(\phi(g)g^{-1})^{-1} -  \phi(\psi(g)g^{-1})(\psi(g)g^{-1})^{-1} \;\; \mbox{mod}\; G_{i+j+2}.
\end{array}
$$
{\it Proof.} By a straightforward computation, we obtain
$$  \begin{array}{l}
[\psi,\phi](g) g^{-1}\\
= [\psi,\phi]((\phi(g)g^{-1})^{-1}) \cdot (\psi \phi \psi^{-1})((\psi(g)g^{-1})^{-1})\cdot \psi(\phi(g)g^{-1}) \cdot \psi(g)g^{-1}.
\end{array}
$$
Since $[\psi,\phi] \in {\rm A}_G(i+j)$ by Proposition 2.2.5 (1) and $\phi(g) g^{-1} \in G_{j+1}$ by Lemma 2.2.4, we have 
$$ [\psi,\phi]((\phi(g)g^{-1})^{-1})  \equiv (\phi(g)g^{-1})^{-1} \;\; {\rm mod} \; G_{i+2j+1}. $$
Similarly, we have
$$ (\psi \phi \psi^{-1})((\psi(g)g^{-1})^{-1}) \equiv \phi((\psi(g) g^{-1})^{-1}) \;\; {\rm mod} \; G_{2i+j+1}. $$
By these three equations together, we have
$$  \begin{array}{l}
[\psi,\phi](g) g^{-1}\\
\equiv  (\phi(g)g^{-1})^{-1} \cdot \phi((\psi(g)(g^{-1})^{-1})\cdot \psi(\phi(g)g^{-1}) \cdot \psi(g)g^{-1} \;\; {\rm mod} \; G_{i+j+2}.
\end{array}
$$
Since $\psi(g) g^{-1} \in G_{i+1}, \phi(g)g^{-1} \in G_{j+1}$ and $[G_{i+1},G_{j+1}]\subset G_{i+j+2}$, we have
$$  \begin{array}{l}
[\psi,\phi](g) g^{-1}\\
\equiv  (\phi(g)g^{-1})^{-1} \cdot \psi(\phi(g)g^{-1}) \cdot \phi((\psi(g) g^{-1})^{-1})\cdot \psi(g)g^{-1} \;\; {\rm mod} \; G_{i+j+2}.
\end{array}
$$ 
Since we easily see that 
$$\left\{ 
\begin{array}{l}
(\phi(g)g^{-1})^{-1} \psi(\phi(g)g^{-1}) \equiv \psi(\phi(g)g^{-1}) (\phi(g)g^{-1})^{-1} \; {\rm mod} \; G_{i+j+2},\\
\phi((\psi(g) g^{-1})^{-1})\cdot \psi(g)g^{-1} \equiv (\phi(\psi(g) g^{-1})\cdot (\psi(g)g^{-1})^{-1})^{-1} \; {\rm mod}\; G_{i+j+2},
\end{array}
\right.
$$
we obtain the assertion. $\;\; \Box$
\\
\\
By Proposition 2.2.5, we can form the graded Lie algebra over $\mathbb{F}_p$ associated to the Andreadakis-Johnson filtration:
$$ {\rm gr}({\rm A}_G) := \bigoplus_{m\geq 0} {\rm gr}_m({\rm A}_G), \;\; {\rm gr}_m({\rm A}_G) := {\rm A}_G(m)/{\rm A}_G(m+1), $$
where the Lie bracket is given by the commutator on the group ${\rm Aut}(G)$. \\
Then by Lemma 2.2.11, the direct sum of Johnson homomorphisms $\tau_m$ over all $m \geq 1$ defines a Lie algebra homomorphism from ${\rm gr}({\rm A}_G)$ to the derivation algebra of ${\rm gr}(G)$ as follows. Recall that an $\mathbb{F}_p$-linear endomorphism of ${\rm gr}(G)$ is called a {\it derivation} on ${\rm gr}(G)$ if it satisfies 
$$\delta([x,y]) = [\delta(x),y]+[x,\delta(y)] \;\;\; (x,y \in {\rm gr}(G)).$$
 Let ${\rm Der}({\rm gr}(G))$ denote the associative $\mathbb{F}_p$-algebra of all derivations on ${\rm gr}(G)$ which has  a Lie algebra structure over $\mathbb{F}_p$ with the Lie bracket defined by $[\delta, \delta'] := \delta \circ \delta' - \delta' \circ \delta$ for $\delta, \delta' \in {\rm Der}({\rm gr}(G))$.  For $m \geq 0$, we define the subspace ${\rm Der}_m({\rm gr}(G))$ of ${\rm Der}({\rm gr}(G))$, the degree $m$ part,  by
$$ {\rm Der}_m({\rm gr}(G)):= \{ \delta \in {\rm Der}({\rm gr}(G)) \, | \, \delta({\rm gr}_n(G)) \subset {\rm gr}_{m+n}(G) \; {\rm for}\; n \geq 1 \}$$
so that we have
$$ {\rm Der}({\rm gr}(G)) = \bigoplus_{m\geq 0} {\rm Der}_m({\rm gr}(G)).$$
Since a derivation on ${\rm gr}(G)$ is determined by its restriction on $H = {\rm gr}_1(G)$, we have a natural inclusion
$$ {\rm Der}_m({\rm gr}(G)) \subset  {\rm Hom}_{\mathbb{F}_p}(H,{\rm gr}_{m+1}(G)); \;\; \delta \mapsto \delta|_H $$
for each $m \geq 1$ and hence we have the inclusion
$$ {\rm Der}_+({\rm gr}(G)) \subset \bigoplus_{m \geq 1} {\rm Hom}_{\mathbb{F}_p}(H, {\rm gr}_{m+1}(G)),$$
where ${\rm Der}_+({\rm gr}(G))$ is the Lie subalgebra of ${\rm Der}({\rm gr}(G))$ consisting of positive degree parts.\\
\\
{\bf Proposition 2.2.12.} {\it The direct sum of $\tau_m$ over $m\geq 1$ defines the Lie algebra homomorphism}
$$ {\rm gr}(\tau) := \bigoplus_{m\geq 1} \tau_m \, : \, {\rm gr}({\rm A}_G) \longrightarrow {\rm Der}_+({\rm gr}(G)).$$
{\it Proof.} (cf. [Da; Proposition 3.18]) By Lemma 2.2.11, it suffices to show that for $\phi \in {\rm A}_G(m)$, the map $g \mapsto \phi(g)g^{-1}$ is indeed a derivation on ${\rm gr}(G)$.  Let  $\phi \in {\rm A}_G(m)$ $(m \geq 1)$ and $g \in G_i, h \in G_j$. By using the commutator formulas
$$ [ab,c] = a[b,c]a^{-1}\cdot[a,c], \;\; [a,bc] = [a,b]\cdot  b[a,c]b^{-1} \;\; (a,b,c \in G),$$
we obtain
$$ \begin{array}{l}
\phi([g,h]) [g,h]^{-1} \\
 = [\phi(g),\phi(h)][g,h]^{-1}\\
 = [gg^{-1}\phi(g),\phi(h)h^{-1}h][g,h]^{-1}\\
 =  g( [g^{-1}\phi(g),\phi(h)h^{-1}]\cdot (\phi(h)h^{-1})[ g^{-1}\phi(g),h](\phi(h)h^{-1})^{-1})g^{-1}\\
  \;\; \;\;\; \cdot [g,\phi(h)h^{-1}] (\phi(h)h^{-1})[g,h](\phi(h)h^{-1})^{-1}[g,h]^{-1}\\
 = g( [g^{-1}\phi(g),\phi(h)h^{-1}]\cdot (\phi(h)h^{-1})[ g^{-1}\phi(g),h](\phi(h)h^{-1})^{-1})g^{-1}\\
 \;\; \;\;\; \cdot [g,\phi(h)h^{-1}] [\phi(h)h^{-1},[g,h]]. 
\end{array}
$$
Since $g^{-1}\phi(g) \in G_{i+m}, \phi(h)h^{-1} \in G_{j+m}$ by Lemma 2.2.4, we have
$$ [g^{-1}\phi(g), \phi(h)h^{-1}] \in G_{i+j+2m}. $$
Similarly, we have
$$ [\phi(h)h^{-1},[g,h]] \in G_{i+2j+m}. $$
By these three claims together, we have
$$ \begin{array}{l}
\phi([g,h]) [g,h]^{-1} \\
\equiv g \phi(h)h^{-1} [ g^{-1}\phi(g),h](g \phi(h)h^{-1})^{-1}
 [g,\phi(h)h^{-1}] \;\; {\rm mod}\; G_{i+j+m+1}.
\end{array}$$
Noting $x [g^{-1}\phi(g),h] x^{-1} \equiv [ g^{-1}\phi(g),h] \; {\rm mod} \; G_{i+j+m+1}$ for $x \in G$, our claim is proved.
$\;\; \Box$
\\

\begin{center}
{\bf 3. Non-Abelian Iwasawa theory}
\end{center}

In this section, we propose an approach to  study non-Abelian Iwasawa theory by means of the Johnson homomorphisms. In the course, we introduce some  invariants from a dynamical viewpoint. 

Throughout this section,  a fixed prime number $p$ is assumed to be odd.
\\
\\
{\bf 3.1. Classical Iwasawa theory.} Let $k$ be a number field of finite degree over $\mathbb{Q}$ and let $k_{\infty}$ be a $\mathbb{Z}_p$-extension of $k$, namely, $k_{\infty}/k$ is a Galois extension whose Galois group is isomorphic to the additive group of $p$-adic integers $\mathbb{Z}_p$.  We call $k_\infty$ the {\it cyclotomic} $\mathbb{Z}_p$-extension of $k$  if $k_\infty$ is the unique $\mathbb{Z}_p$-extension of $k$ contained in $k(\mu_{p^\infty})$.  Let $S_p$ denote the set of primes of $k$ lying over $p$ and $S$ a finite set of primes of $k$ containing $S_p$. Note that the extension $k_\infty/k$ is unramified outside $S_p$. Let $\tilde{k}_S$ be the maximal pro-$p$ extension of $k$ which is unramified outside $S$, and let $M$ be a subextension of $\tilde{k}_S/k$ such that $M/k$ is a Galois extension. We set 
$$ \Gamma := {\rm Gal}(k_\infty/k),\; {\cal G} := {\rm Gal}(M/k) \; \mbox{and}\;  G := {\rm Gal}(M/k_\infty) \leqno{(3.1.1)}$$
so that we have the exact sequence
$$ 1 \longrightarrow G \longrightarrow {\cal G} \longrightarrow \Gamma \longrightarrow 1. \leqno{(3.1.2)}$$
We assume that $G$ is a finitely generated pro-$p$ group, in other words, the $\mu$-invariant is zero.\\

We fix a topological generator $\gamma$ of $\Gamma$ and its lift $\tilde{\gamma} \in {\cal G}$. We then define the automorphism $\phi_{\tilde{\gamma}}$ of $G$ by ${\rm Inn}(\tilde{\gamma})$
$$ \phi_{\tilde{\gamma}}(g) = \tilde{\gamma}g\tilde{\gamma}^{-1} \;\; (g \in G). \leqno{(3.1.3)}$$
We note that if we choose a different lift $\tilde{\gamma}'$ of $\gamma$, $\phi_{\tilde{\gamma}'}$ differs from $\phi_{\tilde{\gamma}}$ by an inner automorphism of $G$ :
$$ \phi_{\tilde{\gamma}'} = {\rm Inn}(x) \circ \phi_{\tilde{\gamma}} \;\; (x := \tilde{\gamma}' \tilde{\gamma}^{-1} \in G). \leqno{(3.1.4)}$$

Let $H$ be the Frattini quotient of $G$, $H = G/G^p[G,G]$, as in (2.1.5). The $\mathbb{F}_p$-linear automorphism $[\phi_{\tilde{\gamma}}]$ of $H$ induced by $\phi_{\tilde{\gamma}}$ is independent of the choice of a lift $\tilde{\gamma}$ and so is denoted by $[\phi_{\gamma}]$. Similarly, we let $H_\infty$ be the Abelianization of $G$,  $H_\infty = G/[G,G]$,   and $[\phi_{\gamma}]_\infty$  the $\mathbb{Z}_p$-module automorphism of $H_\infty$ induced by $\phi_{\tilde{\gamma}}$, which  is independent of the choice of a lift $\tilde{\gamma}$ of $\gamma$. The reason that we use the Zassenhaus filtrarion instead of the lower central series throughout this paper is that any $p$-power of  $\phi_{\tilde{\gamma}}$ acts non-trivially on  $G/[G,G]$ in general. \\

By the Magnus correspondence $\gamma \mapsto 1+X$, we identify the complete group algebra $\mathbb{F}_p[[\Gamma]]$ (resp. $\mathbb{Z}_p[[\Gamma]]$) with the power series algebra $\mathbb{F}_p[[X]]$ (resp. $\mathbb{Z}_p[[X]]$). We set simply $\Lambda := \mathbb{Z}_p[[X]]$ (Iwasawa algebra) and $\overline{\Lambda} := \mathbb{F}_p[[X]]$. Classical Iwasawa theory studies the $\Lambda$-module structure of $H_\infty$, in other words, the $p$-power iterated action of $[\phi_\gamma]_\infty$ on $H_\infty$.
A fundamental theorem of Iwasawa ([Iw1]), under our assumption on $G$, tells us that there is a $\Lambda$-module homomorphism, called a {\it pseudo-isomorphism}, 
$$ H_\infty \longrightarrow \bigoplus_{i=1}^s \Lambda/(f_i(X))  \leqno{(3.1.5)}$$
with finite kernel and cokernel, where $f_i(X)$ is a power of an irreducible distinguished polynomial. Recall that a nonconstant polynomial $f(X) \in \mathbb{Z}_p[X]$ is called {\it distinguished} if $f(X)$ has the form $X^d + a_1X^{d-1}+\cdots + a_d$ with all $a_i \equiv 0$ mod $p$. The {\it Iwasawa polynomial} ($p$-adic zeta function) associated to $H_\infty$ is defined by $\prod_{i=1}^s f_i(X)$.   The set of degrees of $f_i$, $\{ {\rm deg}(f_1), \dots , {\rm deg}(f_s) \}$,   is also an invariant of the $\Lambda$-module $H_\infty$. The Iwasawa $\lambda$-invariant $\lambda(H_\infty)$ is defined by their sum $\sum_{i=1}^s {\rm deg}(f_i)$. 

In some cases,  the pseudo-isomorphism in (3.1.5) turns out to be an isomorphism. Then we can describe the $p$-power iterated action of $[\phi_{\gamma}]$ on $H$ in terms of ${\rm deg}(f_i)$'s. 
Since $H$ is finite, there is an integer $d \geq 0$ such that $[\phi_{\gamma}]^{p^d} = [\phi_{\gamma^{p^d}}] = {\rm id}_H$, namely, $[\phi_\gamma]^{p^d} \in {\rm IA}(G)$. We call such smallest integer $d$ the {\it $p$-period} of $[\phi_\gamma]$ on $H$.
\\
\\
{\bf Proposition 3.1.6.} {\it Suppose that we have a $\Lambda$-module isomorphism
$$ H_\infty \simeq \bigoplus_{i=1}^s \Lambda/(f_i(X)),$$
where $f_i$ is a distinguished polynomial of degree ${\rm deg}(f_i)$. Let $d(H_\infty)$ denote the maximum of ${\rm deg}(f_1), \dots , {\rm deg}(f_s)$. Then we have
$$ [\phi_\gamma]^{p^d} =  [\phi_{\gamma^{p^d}}] = {\rm id}_H, \; \mbox{namely},  [\phi_\gamma]^{p^d} \in {\rm IA}(G)$$
if and only if 
$$ p^d \geq  d(H_\infty).$$
Hence the $p$-period of $[\phi_{\gamma}]$ is  given by the smallest integer $\geq  \log_p d(H_\infty).$}\\
\\
{\it Proof.}  By the assumption, we have a $\overline{\Lambda}$-module isomorphism
$$ H \simeq \bigoplus_{i=1}^s \overline{\Lambda}/(X^{{\rm deg}(f_i)}).$$
Since the action of $[\phi_{\gamma}]^{p^d} - {\rm id}_H$ on $H$ corresponds the multiplication  by $(1+X)^{p^d}-1 = X^{p^d}$, 
$[\phi_{\gamma}]^{p^d} = {\rm id}_H$ if and only if $X^{p^d} \in (X^{{\rm deg}(f_i)})$ for all $i$. From this the assertion follows. $\;\; \Box$
\\
\renewcommand{\thefootnote}{\fnsymbol{footnote}}
\\
{\bf Example 3.1.7.}\footnote[1]{We thank Y. Mizusawa for informing us of this example.}  Let $k := \mathbb{Q}(\mu_p), k_{\infty} := \mathbb{Q}(\mu_{p^\infty})$ and $M$ the maximal unramified pro-$p$ extension of $k_\infty$.  The assumption of Proposition 3.1.5 is then satisfied if the Vandiver conjecture is true, namely, $p$ does not divide the class number of the maximal real subfield of $k$ ([Wa; Theorem 10.16]). The Vandiver conjecture is known to be true for $p < 163 577 856$ ([BH]). For instance, we have $H_\infty = \Lambda/(f)$ for $p = 37$ and $H_\infty = \Lambda/(f_1) \oplus \Lambda/(f_2)$ for $p = 157$, where $f, f_1$ and $f_2$ are all distinguished polynomials of degree one ([IS]). So, the $p$-period of $[\phi_\gamma]$ is zero, namely,  $[\phi_\gamma]$ acts trivially on $H$. Mizusawa made a program to compute the Iwasawa polynomial when $k$ is an imaginary quadratic field $\mathbb{Q}(\sqrt{-D})$, $k_\infty$ is the cyclotomic $\mathbb{Z}_p$-extension and $M$ is the maximal unramified pro-$p$ extension of $k_\infty$. For example, when $p = 3$ and  $D = 186, 211, 231, 249$, $H_\infty = \Lambda/(f)$ with ${\rm deg}(f) = 2$ and so the $3$-period  of $[\phi_\gamma]$ is one, and when $p=3$ and $D = 214, 274$, $ H_\infty = \Lambda/(f)$ with ${\rm deg}(f) = 4$ and so the $3$-period of $[\phi_\gamma]$ is two.\\
\\
{\bf  3.2. Non-Abelian Iwasawa theory via Johnson homomorphisms.}  A basic problem in non-Abelian Iwasawa theory is to understand the $p$-power iterated action of $\phi_{\tilde{\gamma}}$ on $G$, while classical Iwasawa theory deals with that of $[\phi_\gamma]$ on $H_{\infty}$ as shown in 3.1. Let $\{ G_n \}_{n\geq 1}$ be the Zassenhaus filtration of $G$ so that $H = G/G_2$, and let $[\phi_{\tilde{\gamma}}]_m$ be the automorphism of $G/G_{m+1}$ induced by $\phi_{\tilde{\gamma}}$ as defined in (2.2.1). We aim to study the $p$-power iterated action of $[\phi_{\tilde{\gamma}}]_m$ on $G/G_{m+1}$ for all $m \geq 1$ by means of the $p$-Johnson homomorphisms introduced in 2.2. \\

First, let us see how a different choice of a lift of $\gamma$ affects the  action of a power of $[\phi_{\tilde{\gamma}}]_m$ on $G/G_{m+1}$ \\
\\
{\bf Lemma 3.2.1.} {\it Let $\tilde{\gamma}$, $\tilde{\gamma}'$ be lifts of $\gamma$ in ${\cal G}$ and set $x = \tilde{\gamma}' \tilde{\gamma}^{-1} \in G$ as in $(3.1.4)$. Suppose $x \in G_m$. Then, for each integer $e \geq 1$, we have}
$$ \phi_{\tilde{\gamma}'}^e \in {\rm A}_G(m) \Longleftrightarrow  \phi_{\tilde{\gamma}}^e \in {\rm A}_G(m).$$
{\it Proof.}  By (3.1.4), we have
$$ \phi_{\tilde{\gamma}'}^e(g) = y \phi_{\tilde{\gamma}}^e(g) y^{-1}, \;\; y := x \phi_{\tilde{\gamma}}(x)\cdots \phi_{\tilde{\gamma}}^{e-1}(x) \in G_m,$$
for any $g \in G$. Since elements of $G/G_{m+1}$ and $G_m/G_{m+1}$ commute, the assertion is shown as follows
$$ \begin{array}{ll}
\phi_{\tilde{\gamma}'}^e \in {\rm A}_G(m) &  \Leftrightarrow \phi_{\tilde{\gamma}'}^e(g) g^{-1} \in G_{m+1} \; \mbox{for any}\; g \in G\\
                                                                     & \Leftrightarrow y \phi_{\tilde{\gamma}}^e(g) y^{-1} g^{-1} \in G_{m+1} \; \mbox{for any}\; g \in G\\
                                                                    & \Leftrightarrow  \phi_{\tilde{\gamma}}^e(g) g^{-1} \in G_{m+1} \; \mbox{for any}\; g \in G\\
                                                                    & \Leftrightarrow \phi_{\tilde{\gamma}} \in {\rm A}_G(m) \;\; \Box
\end{array}
$$
\\
 Let  ${\rm gr}(G) = \bigoplus_{n\geq 1} {\rm gr}_n(G)$, ${\rm gr}_n(G) = G_n/G_{n+1}$, be the graded  Lie algebra over $\mathbb{F}_p$ associated to the Zassenhaus filtration of $G$ as in (2.1.6), and let $\{{\rm A}_G(m)\}_{m \geq 0}$ be the Andreadakis-Johnson filtration of ${\rm Aut}(G)$. For $m \geq 1$, let
 $$ \tau_m \, : \, {\rm A}_G(m) \longrightarrow {\rm Hom}_{\mathbb{F}_p}(H, {\rm gr}_{m+1}(G))$$
 be the $p$-Johnson homomorphism. The next Corollary follows immediately from Lemma 3.2.1.
\\
\\
{\bf Corollary 3.2.2.} {\it Let $\tilde{\gamma}$, $\tilde{\gamma}'$ be lifts of $\gamma$ in ${\cal G}$ and set $x = \tilde{\gamma}' \tilde{\gamma}^{-1} \in G$. Suppose $x \in G_{m+1}$ and $\phi_{\tilde{\gamma}}^e \in {\rm A}_G(m)$ $(e \geq 1)$.
Then we have}
$$ \tau_m(\phi_{\tilde{\gamma}'}^e) = \tau_m (\phi_{\tilde{\gamma}}^e).$$
{\it Proof.} By Lemma 3.2.1, $\phi_{\tilde{\gamma}'}^e \in {\rm A}_G(m)$. Since $ \phi_{\tilde{\gamma}'}^e = {\rm Inn}(y) \circ \phi_{\tilde{\gamma}}^e$ with $y =  x \phi_{\tilde{\gamma}}(x)\cdots \phi_{\tilde{\gamma}}^{e-1}(x) \in G_{m+1}$, the assertion follows from Theorem 2.2.8 and Proposition 2.2.10. $\;\; \Box$\\

We fix a lift $\tilde{\gamma} \in {\cal G}$ of $\gamma$. Generalizing the $p$-period of $[\phi_{\gamma}]$ on $H = G/G_m$, we define  the {\it $p$-period} $d(m)$ of $\phi_{\tilde{\gamma}}$ acting on $G/G_{m+1}$ for each $m \geq 1$  by the smallest integer $d \geq 0$ such that 
$$ \phi_{\tilde{\gamma}}^{p^{d}}  \in  {\rm A}_G(m). \leqno{(3.2.3)}$$
Thus we have non-decreasing sequence $\{ d(m)\}_{m\geq 1}$ of integers. \\
\\
{\bf Lemma 3.2.4.} {\it For each integer $m \geq 1$, we have}
$$ d(m+1) = d(m) \; {\rm or}\; d(m) +1.$$
{\it Proof.}  By definition of $d(m)$, we have $d(m+1) \geq d(m)$. Suppose $\phi_{\tilde{\gamma}}^{p^{d}}  \in  {\rm A}_G(m)$. Then by Proposition 2.2.5 (2), we have $\phi_{\tilde{\gamma}}^{p^{d+1}}  \in  {\rm A}_G(m+1)$.
Hence $d(m+1) \leq d(m) + 1$. $\;\; \Box$\\
\\
Now we introduce another sequence of integers $\{ m(d)\}_{d \geq 0}$ as follows. For each integer $d \geq 0$, we define the integer $m(d) \geq 1$ by
$$ \phi_{\tilde{\gamma}}^{p^d} \in {\rm A}_G(m(d)), \;\; \phi_{\tilde{\gamma}}^{p^d} \notin {\rm A}_G(m(d)+1). \leqno{(3.2.5)}$$
It is a strictly increasing sequence. In fact, we have\\
\\
{\bf Lemma 3.2.6.} {\it For each integer $d \geq 0$, we have}
$$ m(d+1) \geq m(d) + 1. $$
{\it Proof.}  Since   $\phi_{\tilde{\gamma}}^{p^{d}}  \in  {\rm A}_G(m(d))$ for each $d \geq 0$, by Proposition 2.2.5 (2), we have  $\phi_{\tilde{\gamma}}^{p^{d+1}}  \in  {\rm A}_G(m(d)+1)$. Hence, by definition (3.2.5), 
we have $m(d+1) \geq m(d)+1$. $\;\; \Box$
\\
\\
Then the sequence $\{ \tau_{m(d)}(\phi_{\tilde{\gamma}}^{p^d})\}_{d \geq 0}$ in ${\rm Hom}_{\mathbb{F}_p}(H,{\rm gr}_{m(d)+1}(G))$ describes the action of $\phi_{\tilde{\gamma}}^{p^d}$ on $G/G_{m(d)+1}$ for all $d \geq 0$. 
 In Section 5, we give a cohomological interpretation of  $\tau_{m(d)}(\phi_{\tilde{\gamma}}^{p^d})$ in terms of Massey products  in Galois cohomology.\\
 \\
 {\bf Remark 3.2.7.} Let $M$ the maximal unramified pro-$p$ extension of $k_\infty$. Ozaki ([O]) studied the $\Gamma$-action on the graded pieces associated to the lower central series of $G = {\rm Gal}(M/k_{\infty})$ and obtained arithmetic results.  We also refer to Sharifi's paper [Sh] for a related work. Our approach is different from theirs.\\

\begin{center}
{\bf 4.  $p$-Johnson maps for a free pro-$p$ group}
\end{center}

In this section, following Kawazumi ([Kw]), we extend the $p$-Johnson homomorphisms in Section 2 to maps defined on the whole group of automorphisms when $G$ is a free pro-$p$ group. 

Throughout this section, let $F$ denote a free pro-$p$ group on $x_1,\dots , x_r$. 
A fixed prime number $p$ is arbitrary in 4.1 and assumed to be odd in 4.2.
\\
\\
{\bf 4.1. $p$-Johnson maps.}  We keep the same notations as in 2.1, only replacing $G$ by $F$. Let $\mathbb{F}_p[[F]]$ be the complete group algebra of $F$ over $\mathbb{F}_p$ with augmentation ideal $I_F$. Let
$\{ F_n \}_{n\geq 1}$ be the Zassenhaus filtration defined by $F_n = F \cap (1+I_F^n)$ and let $H := F/F_2 = F/F^p[F,F]$ be the Frattini quotient of $F$. We write $[f]$ for the image of $f \in F$ in $H$: $[f] := f$ mod $F_2$. We denote  $[x_j]$  by $X_j$ ($1\leq j \leq r$)  simply so that $H$ is a vector space over $\mathbb{F}_p$ with basis $X_1,\dots,X_r$
$$ H = \mathbb{F}_p X_1 \oplus \cdots \oplus  \mathbb{F}_p X_r.$$

As in 2.1, let ${\rm gr}(F)$ be the graded restricted Lie algebra over $\mathbb{F}_p$ associated to the Zassenhaus filtration $\{ F_n \}_{n\geq 1}$ of $F$
$$  {\rm gr}(F) :=  \bigoplus_{n\geq 1} {\rm gr}_n(F),\;\; {\rm gr}_n(F) := F_n/F_{n+1}.$$
It is the free Lie algebra over $\mathbb{F}_p$ on $X_1,\dots, X_r$. Its restricted universal enveloping algebra $U {\rm gr}(F)$  is given by the graded associative algebra ${\rm gr}(\mathbb{F}_p[[F]])$ (cf. (2.1.7))
$$ U {\rm gr}(G) = {\rm gr}(\mathbb{F}_p[[F]]) := \bigoplus_{m \geq 0} {\rm gr}_m(\mathbb{F}_p[[F]]),\;\; {\rm gr}_m(\mathbb{F}_p[[F]]) := I_F^m/I_F^{m+1}$$
together with  the injective restricted Lie algebra homomorphism
$$ {\rm gr}(\theta) = \bigoplus_{m \geq 1} \theta_m \, :\, {\rm gr}(F) \longrightarrow  {\rm gr}(\mathbb{F}_p[[F]]),$$
where $ \theta_m \, : \, {\rm gr}_m(F) \rightarrow {\rm gr}_m(\mathbb{F}_p[[F]])$ is given by
$$\theta_m(f\; {\rm mod}\; F_{m+1}) := f -1 \; {\rm mod}\; I_F^{m+1}.$$
By the correspondence $ x_j - 1 \; {\rm mod}\; I_F^2  \in {\rm gr}_1(\mathbb{F}_p[[F]]) \mapsto X_j \in H$, the universal envelope ${\rm gr}(\mathbb{F}_p[[F]])$ is identified with the tensor algebra on $H$ over $\mathbb{F}_p$ or the non-commutative polynomial algebra $\mathbb{F}_p\langle X_1,\dots ,X_r \rangle$ of variables $X_1,\dots ,X_r$ over $\mathbb{F}_p$
$$ \begin{array}{ll}
U {\rm gr}(G) = {\rm gr}(\mathbb{F}_p[[F]]) & = \, \displaystyle{\bigoplus_{m \geq 0} H^{\otimes m}}\\
                                               &= \, \mathbb{F}_p\langle X_1,\dots ,X_r \rangle.
\end{array}
$$   
Here the graded piece ${\rm gr}_m(\mathbb{F}_p[[F]])$ corresponds to $H^{\otimes m}$, the vector space over $\mathbb{F}_p$ with basis $X_{i_1}\cdots X_{i_m}$ $(1\leq i_1,\dots ,i_m \leq r)$, monomials  of degree $m$,
and so $\theta_m$ may be regarded as the injective $\mathbb{F}_p$-linear map 
$$ \theta_m \, : \, {\rm gr}_m(F) \hookrightarrow H^{\otimes m} \leqno{(4.1.1)}$$

In order to extend the Johnson homomorphisms in 2.2 to the maps defined on the whole automorphism group ${\rm Aut}(F)$, we work with the completion $\widehat{U}$ of the universal envelope $ U {\rm gr}(F) = {\rm gr}(\mathbb{F}_p[[F]])$ with respect to $I_F$-adic topology. So $\widehat{U}$ is the complete tensor algebra on $H$ over $\mathbb{F}_p$ which is identified with the $\mathbb{F}_p$-algebra $\mathbb{F}_p\langle \langle X_1,\dots ,X_r \rangle\rangle$ of non-commutative formal power series of  variables $X_1,\dots,X_r$ over $\mathbb{F}_p$ (In [Kw] Kawazumi wrote $\widehat{T}$ for $\widehat{U}$)
$$ \begin{array}{ll}       \widehat{U}  & := \displaystyle{\prod_{m \geq 0} H^{\otimes m}}\\
                                                          &    = \mathbb{F}_p\langle \langle X_1,\dots ,X_r \rangle\rangle.
\end{array}$$             
Then the composite of $\theta_m$  in (4.1.1) with the natural inclusion $H^{\otimes m} \hookrightarrow \widehat{U}$ is nothing but the restriction to $F_m$ of the {\it Magnus embedding}
$$ \theta \, : \, F \hookrightarrow  \widehat{U}^{\times} \leqno{(4.1.2)}$$
defined by $\theta(x_j) := 1+X_j$ $(1\leq j \leq r)$. \\

For $n \geq 1$, we let
 $$\widehat{U}_n := \prod_{m \geq n} H^{\otimes m}$$ 
be the two-sided ideal of $\widehat{U}$ corresponding to formal power series of degree $\geq n$.  An $\mathbb{F}_p$-algebra  automorphism $\varphi$ of $\widehat{U}$ is then called {\it filtration-preserving} if $\varphi(\widehat{U}_n) = \widehat{U}_n$ for all $n \geq 0$ and we denote by ${\rm Aut}^{\rm fil}(\widehat{U})$ the group of filtration-preserving $\mathbb{F}_p$-algebra automorphisms of $\widehat{U}$. The following useful Lemma, which we call {\it Kawazumi's lemma},  gives a criterion for a $\mathbb{F}_p$-algebra endomorphism of $\widehat{U}$ to be a filtration-preserving automorphism.
\\
\\
{\bf Lemma 4.1.3.}  ({\it Kawazumi's lemma}).  {\it A $\mathbb{F}_p$-algebra endomorphism $\varphi$ of $\widehat{U}$ is a  filtration-preserving automorphism of $\widehat{U}$, $\varphi \in {\rm Aut}^{\rm fil}(\widehat{U})$, if and only if the following conditions are satisfied}:
\\
(1)  {\it $\varphi(\widehat{U}_n) \subset \widehat{U}_n$ for all $n \geq 0$.}\\
(2)  {\it the induced $\mathbb{F}_p$-linear map $[\varphi]$ on $\widehat{U}_1/\widehat{U}_2 = H$ defined by $[\varphi](h) := \varphi(h) \, {\it mod}\, \widehat{U}_2$ $(h \in H)$ is an isomorphism.}\\
\\
{\it Proof.}  Suppose $\varphi \in {\rm Aut}^{\rm fil}(\widehat{U})$. Since $\varphi$ is filtration-preserving, the condition (1) holds. To show the condition (2), consider the following commutative diagram for vector spaces over $\mathbb{F}_p$ with exact rows:
$$ \begin{array}{ccccccccc}
0 & \longrightarrow & \widehat{U}_2 & \longrightarrow & \widehat{U}_1 & \longrightarrow & H & \longrightarrow & 0 \\
  &                          & \;\;\;\;\;\;\; \downarrow {\footnotesize \varphi|_{\widehat{U}_2}} & & \;\;\;\;\;\;\; \downarrow  {\footnotesize \varphi|_{\widehat{U}_1}} &  & \;\;\;\;\; \downarrow {\footnotesize [\varphi]} & & \\
  0 & \longrightarrow & \widehat{U}_2 & \longrightarrow & \widehat{U}_1 & \longrightarrow & H & \longrightarrow & 0. 
  \end{array}
  $$
  Since $\varphi(\widehat{U}_n) = \widehat{U}_n$ for all $n \geq 0$, we have ${\rm Coker}(\varphi|_{\widehat{U}_i}) = 0$ for $i=1,2$, in particular. Since $\varphi$ is an automorphism, we have ${\rm Ker}(\varphi) = 0$, in particular, 
  ${\rm Ker}(\varphi|_{\widehat{U}_i}) = 0$ for $i=1,2$. By snake lemma applied to the above diagram, we obtain ${\rm Ker}([\varphi]) = 0$ and ${\rm Coker}([\varphi]) = 0$, hence the condition (2).

Suppose  that  an $\mathbb{F}_p$-algebra endomorphism $\varphi$ of $\widehat{U}$ satisfies the conditions (1) and (2). Let $z = (z_m)$ be any element of $\widehat{U}$ with $z_m \in H^{\otimes m}$ for $m \geq 0$. To show that $\varphi$ is an automorphism, we have only to prove that there exists uniquely $y = (y_m) \in \widehat{U}$ such that
$$ z = \varphi(y). \leqno{(4.1.3.1)}$$
 Note by the condition (1) and (2) that $\varphi$ induces an $\mathbb{F}_p$-linear automorphism of $\widehat{U}_m/\widehat{U}_{m+1} = H^{\otimes m}$, which is nothing but $[\varphi]^{\otimes m}$. Then, writing $\varphi(y_i)_j$ for the component
 of $\varphi(y_i)$ in $H^{\otimes j}$ for $i < j$, the equation (4.1.3.1) is equivalent to the following system of equations:
 $$\left\{ 
\begin{array}{l}
  z_0 = \varphi(y_0) = y_0,\\
 z_1 = [\varphi](y_1),\\
 z_2 = [\varphi]^{\otimes 2}(y_2) + \varphi(y_1)_2,\\
 \cdots \\
 z_m = [\varphi]^{\otimes m}(y_m) + \varphi(y_1)_m + \cdots + \varphi(y_{m-1})_m ,\\
 \cdots 
 \end{array} 
\right.
  \leqno{(4.1.3.2)}
 $$
Since $[\varphi]^{\otimes m}$ is an automorphism, we can find the unique solution $y = (y_m)$ of  (4.1.3.2)  from the lower degree. Therefore $\varphi$ is an $\mathbb{F}_p$-algebra automorphism. Furthermore, we can see easily that if $z_0=\cdots = z_{n-1}=0$, then $y_0=\cdots = y_{n-1}=0$ for $n \geq 1$. This means that $\varphi^{-1}(\widehat{U}_n) \subset \widehat{U}_n$ and so $\varphi$ is filtration-preserving. $\;\; \Box$
\\

By Lemma 4.1.3, each $\varphi \in {\rm Aut}^{\rm fil}(\widehat{U})$ induces an $\mathbb{F}_p$-linear automorphism $[\varphi]$ of $H = \widehat{U}_1/\widehat{U}_2$ and so we have a group homomorphism 
$$ [\;\;] \; : \; {\rm Aut}^{\rm fil}(\widehat{U}) \longrightarrow {\rm GL}(H).$$
We define the {\it induced automorphism group} of $\widehat{U}$ by 
$$ {\rm IA}(\widehat{U}) := {\rm Ker}([\;\;]).$$
We note that there is a natural splitting $s : {\rm GL}(H) \rightarrow {\rm Aut}^{\rm fil}(\widehat{U})$ of $[\;\; ],$ which is defined by 
$$  s(P)((z_m)) := (P^{\otimes m}(z_m)) \;\; \mbox{for} \; P \in {\rm GL}(H).$$
In the following, we also regard $[P] \in {\rm GL}(H)$ as an element of ${\rm Aut}^{\rm fil}(\widehat{U})$ through the splitting $s$ and write simply $[P]$ for $s([P])$. 
Thus we have the following\\
\\
{\bf Lemma 4.1.4.}  {\it We have a semi-direct decomposition}
$$ {\rm Aut}^{\rm fil}(\widehat{U}) \; = \; {\rm IA}(\widehat{U}) \rtimes {\rm GL}(H)$$
{\it given by } $\varphi  =  (\varphi \circ [\varphi]^{-1}, [\varphi])$.\\

Let $\varphi \in {\rm IA}(\widehat{U})$. Since $\varphi$ acts on $\widehat{U}_1/\widehat{U}_2 = H$ trivially, we have
$$ \varphi(h) - h \in \widehat{U}_2 \; \; \mbox{for any} \; h \in H, $$
and so we have a map
$$ E \; : \; {\rm IA}(\widehat{U}) \longrightarrow {\rm Hom}_{\mathbb{F}_p}(H,\widehat{U}_2); \; \varphi \mapsto \varphi|_H - {\rm id}_H,$$
where ${\rm Hom}_{\mathbb{F}_p}(H,\widehat{U}_2)$ denotes the group of $\mathbb{F}_p$-linear maps $H \rightarrow \widehat{U}_2$. The following Proposition will play a key role in our discussion.\\
\\
{\bf Proposition 4.1.5.}  {\it The map $E$ is bijective.}\\
\\
{\it Proof.}  Injectivity: Suppose $E(\varphi) = E(\varphi')$ for $\varphi, \varphi' \in {\rm IA}(\widehat{U})$. Then we have $\varphi|_H = \varphi'|_H$. Since an $\mathbb{F}_p$-algebra endomorphism of $\hat{U}$ is determined by its restriction on $H$, we have $\varphi = \varphi'$. \\
Surjectivity: Take any $\eta \in {\rm Hom}_{\mathbb{F}_p}(H,\widehat{U}_2)$.  We can extend $\eta + {\rm id}_H : H \rightarrow \widehat{U}_2$ uniquely to a $\mathbb{F}_p$-algebra endomorphism $\varphi$ of $\widehat{U}$. Then we have obviously $\varphi(\widehat{U}_n) \subset \widehat{U}_n$ for all $n \geq 0$. Since $\widehat{U}_1/\widehat{U}_2 = H$ and we see that 
$$ [\varphi](h \, {\rm mod} \, \widehat{U}_2) = \varphi(h) \, {\rm mod}\, \widehat{U}_2 = h + \eta(h) \, {\rm mod} \, \widehat{U}_2 = h \, {\rm mod}\, \widehat{U}_2,$$
we have $[\varphi] = {\rm id}_{H}$. By Kawazumi's Lemma 2.1, we have $\varphi \in {\rm IA}(\widehat{U})$ and $E(\varphi) = \eta$.  $\;\; \Box$\\
\\
By  Lemma 4.1.4 and Proposition 4.1.5, we have the following\\
\\
{\bf Corollary 4.1.6.}  {\it We have a bijection}
$$ \hat{E} \; : \; {\rm Aut}^{\rm fil}(\widehat{U}) \simeq {\rm Hom}_{\mathbb{F}_p}(H,\widehat{U}_2) \times {\rm GL}(H)$$
given by $\hat{E}(\varphi) = (E(\varphi \circ [\varphi]^{-1}), [\varphi]).$\\

 The Magnus embedding $\theta : F \hookrightarrow \widehat{U}^{\times}$ in (4.1.2) is extended to an $\mathbb{F}_p$-algebra isomorphism, denoted by the same $\theta$,
$$ \theta : \mathbb{F}_p[[F]] \stackrel{\sim}{\longrightarrow} \widehat{U},  \leqno{(4.1.7)} $$
which satisfies 
$$ \theta(I_F^n) = \widehat{U}_n \;\; \mbox{for}\; m \geq 1.   \leqno{(4.1.8)} $$
For $m \geq 0$, let $\theta_m$ denote the component of $\theta$ in $H^{\otimes m}$ as in (4.1.1):
$$ \theta(\alpha) = \sum_{m=0}^\infty \theta_m(\alpha), \;\; \theta_m(\alpha) \in H^{\otimes m} \;\; (\alpha \in \mathbb{F}_p[[F]]).$$
Note that  $\theta_0(f) = 1$ and $\theta_1(f) = [f]$ for $f \in F$. Further we can write $\theta_m(\alpha)$ as 
$$ \theta_m(\alpha) = \sum_{1\leq i_1,\dots, i_m \leq r} \epsilon(i_1\cdots i_m;\alpha)X_{i_1}\cdots X_{i_m}, \leqno{(4.1.9)}$$
where the coefficient $\epsilon(i_1\cdots i_m;\alpha)$ is given in terms of the pro-$p$ Fox free derivative $\partial /\partial x_j : \mathbb{Z}_p[[F]] \rightarrow \mathbb{Z}_p[[F]]$ ([Ih], [Ms2, 8.3])
$$ \epsilon(i_1\cdots i_m;\alpha) = \epsilon_{\mathbb{Z}_p[[F]]}\left( \frac{\partial^m \tilde{\alpha}}{\partial x_{i_1}\cdots \partial x_{i_m}} \right) \;\; {\rm mod}\; p,
$$
where $\epsilon_{\mathbb{Z}_p[[F]]} : \mathbb{Z}_p[[F]] \rightarrow \mathbb{Z}_p$ is the augmentation map and $\tilde{\alpha} \in \mathbb{Z}_p[[F]]$ such that $\tilde{\alpha} \; {\rm mod}\; p = \alpha$.

An $\mathbb{F}_p$-algebra automorphism $\varphi$ of $\mathbb{F}_p[[F]]$ is said to be {\it filtration-preserving} if $\varphi(I_F^n) = I_F^n$ for all $n \geq 0$ and we denote by ${\rm Aut}^{\rm fil}(\mathbb{F}_p[[F]])$ the group of filtration-preserving automorphisms of $\mathbb{F}_p[[F]]$. By (4.1.7) and (4.1.8), we have an isomorphism
$$ {\rm Aut}^{\rm fil}(\mathbb{F}_p[[F]]) \simeq {\rm Aut}^{\rm fil}(\widehat{U}); \; \varphi \mapsto \theta \circ \varphi \circ \theta^{-1}. \leqno{(4.1.10)}$$

Now, let $\phi \in {\rm Aut}(F)$. Then $\phi$ induces a filtration-preserving $\mathbb{F}_p$-algebra automorphism $\hat{\phi}$ of $\mathbb{F}_p[[F]])$. In fact, 
$\phi$ induces an automorphism $[\phi]_n$ of a finite $p$-group $F/F_n$ and hence an $\mathbb{F}_p$-algebra automorphism, denoted by the same $[\phi]_{n}$,  of a finite group ring $\mathbb{F}_p[F/F_n]$
 $$ [\phi]_{n} \; : \; \mathbb{F}_p[F/F_n] \stackrel{\sim}{\longrightarrow} \mathbb{F}_p[F/F_n]$$
 for each $n \geq 1$, which sends the augmentation ideal of $\mathbb{F}_p[F/F_p]$ onto itself. Taking the inverse limit with respect to $n$, we obtain an $\mathbb{F}_p$-algebra automorphism
 $$ \widehat{\phi} := \lim_{{\scriptstyle \longleftarrow} \atop {\scriptstyle n}} [\phi]_{n}  \; : \; \mathbb{F}_p[[F]] \stackrel{\sim}{\longrightarrow} \mathbb{F}_p[[F]]$$
 such that $\widehat{\phi}(I_F) = I_F$. Thus we have an injective homomorphism
 $$ {\rm Aut}(F) \longrightarrow {\rm Aut}^{\rm fil}(\mathbb{F}_p[[F]]); \; \phi \mapsto \widehat{\phi}.$$
By composing with the isomorphism (4.1.10), we obtain an injective homomorphism
$$ \widehat{\kappa}^{\theta} \; : \; {\rm Aut}(F) \longrightarrow {\rm Aut}^{\rm fil}(\widehat{U}); \; \phi \mapsto \theta\circ \widehat{\phi} \circ \theta^{-1}. $$
\\
{\bf Lemma 4.1.11.} {\it Let $[\phi]$ denote the $\mathbb{F}_p$-linear automorphism of $H$ induced by $\phi \in {\rm Aut}(F)$. Then we have}
$$  [\widehat{\kappa}^{\theta}(\phi)] = [\phi] \; \mbox{in} \;\;  {\rm GL}(H).$$
{\it Proof.} We have, for $X_j \in H$  $(1\leq j \leq r),$
$$ \begin{array}{lll}
            \widehat{\kappa}^{\theta}(\phi)                                             & = & (\theta\circ \widehat{\phi} \circ \theta^{-1})(X_j) \\
                                                                                              & = & (\theta \circ \widehat{\phi} \circ \theta^{-1})(\theta(x_j)-1)\\
                                                                                             & = & (\theta \circ \widehat{\phi})(x_j-1)\\
                                                                                             & = & \theta(\phi(x_j)) -1 \\
                                                                                             & \equiv & [\phi(x_j)] \; \mbox{mod}\; \widehat{U}_2\\
                                                                                             & = & [\phi](X_j) \; \mbox{mod}\; \widehat{U}_2.
 \end{array}
$$
Hence we have $[\widehat{\kappa}^{\theta}(\phi)] = [\phi]$. $\;\; \Box$\\
\\
By Lemma 4.1.11, we have, for $\phi \in {\rm Aut}(F)$,
$$ \widehat{\kappa}^{\theta}(\phi) = (\widehat{\kappa}^{\theta}(\phi)\circ [\phi]^{-1}, [\phi])$$
under the semi-direct decomposition ${\rm Aut}^{\rm fil}(\widehat{U}) = {\rm IA}(\widehat{U}) \rtimes {\rm GL}(H)$ of Lemma 4.1.4. We set
$$ \kappa^{\theta}(\phi) := \widehat{\kappa}^{\theta}(\phi) \circ [\phi]^{-1} = \theta \circ \widehat{\phi} \circ \theta^{-1} \circ [\phi]^{-1} \;\; \;\; (\phi \in {\rm Aut}(F)).  \leqno{(4.1.12)}$$
\vspace{.2cm}

Now, we define the {\it extended  $p$-Johnson map}
$$ \widehat{\tau}^{\theta} \; : \; {\rm Aut}(F) \longrightarrow {\rm Hom}_{\mathbb{F}_p}(H,\widehat{U}_2) \rtimes {\rm GL}(H) $$
by  composing $\widehat{\kappa}^{\theta}$ with $\hat{E}$ of Corollary 4.1.6, and we define the {\it $p$-Johnson map}
$$ \tau^{\theta} \; : \; {\rm Aut}(F) \longrightarrow {\rm Hom}_{\mathbb{F}_p}(H,\widehat{U}_2) 
$$
by the composing $\widehat{\tau}^{\theta}$ with the projection on ${\rm Hom}_{\mathbb{F}_p}(H,\widehat{U}_2)$, namely, for $\phi \in {\rm Aut}(F)$,
$$ \tau^{\theta}(\phi) := E(\kappa^{\theta}(\phi)) = \kappa^{\theta}(\phi)|_H - {\rm id}_H.  \leqno{(4.1.13)}$$
For $m \geq 1$, we define the {\it $m$-th $p$-Johnson map} 
$$ \tau_m^{\theta} \; : \; {\rm Aut}(F) \longrightarrow {\rm Hom}_{\mathbb{F}_p}(H, H^{\otimes (m+1)})$$
by the $m$-th component of $\tau^K$:
$$ \tau^{\theta}(\phi) :=  \sum_{m\geq 1}\tau_{m}^{\theta}(\phi) \;\;\;\;  (\phi \in {\rm Aut}(G)).  \leqno{(4.1.14)}$$

Unlike the $p$-Johnson homomorphisms (Theorem 2.2.8), the $p$-Johnson map $\tau^{\theta} = E \circ \kappa^{\theta} : {\rm Aut}(F) \rightarrow {\rm Hom}(H,\widehat{U}_2)$ is no longer a homomorphism. In fact, we have the following\\
\\
{\bf Proposition 4.1.15.} {\it We have}
$$ \kappa^{\theta}(\phi_1 \circ \phi_2) = \kappa^{\theta}(\phi_1)\circ [\phi_1]\circ \kappa^{\theta}(\phi_2)\circ [\phi_1]^{-1}.$$
\\
{\it Proof.} By (4.1.12), we have
$$ \begin{array}{ll}
\kappa^{\theta}(\phi_1 \phi_2) & = \theta \circ (\widehat{ \phi_1 \phi_2})\circ \theta^{-1} \circ [\phi_1 \phi_2]^{-1}\\
                              & = \theta \circ \widehat{\phi}_1 \circ \widehat{\phi}_2 \circ \theta^{-1} \circ [\phi_2]^{-1} \circ [\phi_1]^{-1}\\
                              & = \theta \circ \widehat{\phi}_1 \circ \theta^{-1} \circ [\phi_1]^{-1} \circ [\phi_1] \circ \theta \circ \widehat{\phi}_2 \circ \theta^{-1} \circ [\phi_2]^{-1} \circ [\phi_1]^{-1}\\
                              & = \kappa^{\theta}(\phi_1) \circ [\phi_1] \circ \kappa^{\theta}(\phi_2) \circ [\phi_1]^{-1}.  \;\;\;\;\;\;\;\;\;\; \Box
\end{array}
$$
\\
Proposition 4.1.15 yields an infinite sequence coboundary relations which Johnson maps $\tau_m^{\theta}$ satisfies. Here we give the formulas for $\tau_1^{\theta}$ and $\tau_2^{\theta}$.\\
\\
{\bf Proposition 4.1.16.}  {\it We have}
$$ \begin{array}{ll}
\tau_1^{\theta}(\phi_1 \phi_2) &= \tau_1^{\theta}(\phi_1) + [\phi_1]^{\otimes 2} \circ \tau_1^{\theta}(\phi_2) \circ [\phi_1]^{-1},\\
\tau_2^{\theta}(\phi_1 \phi_2) & = \tau_2^{\theta}(\phi_1) + (\tau_1^{\theta}(\phi_1)\otimes {\rm id}_H + {\rm id}_H \otimes \tau_1^{\theta}(\phi_1)) \circ [\phi_1]^{\otimes 2} \circ \tau_1^{\theta}(\phi_2) \circ [\phi_1]^{-1}\\
& \;\;\;\;\; \;\;\;\;\; + [\phi_1]^{\otimes 3} \circ \tau_2^{\theta}(\phi_2) \circ [\phi_1]^{-1}.
\end{array}
$$
\\
{\it Proof.} By definition (4.1.14), we have
$$ \tau^{\theta}(\phi_1 \phi_2) = \sum_{m \geq 1} \tau_m^{\theta}(\phi_1 \phi_2). \leqno{(4.1.16.1)}$$
On the other hand, by Proposition 4.1.15 and (4.1.13), we have, for $h \in H$,
$$ \begin{array}{ll}
\tau^{\theta}(\phi_1 \phi_2) & = -h + \kappa^{\theta}(\phi_1 \phi_2)(h)\\
           & = -h + (\kappa^{\theta}(\phi_1)\circ [\phi_1] \circ \kappa^{\theta}(\phi_2) \circ [\phi_1]^{-1})(h)\\
           & = -h + ( \kappa^{\theta}(\phi_1)\circ [\phi_1] \circ ({\rm id}_H + \tau^{\theta}(\phi_2)))([\phi_1]^{-1}(h))\\
           & = \displaystyle{-h + (\kappa^{\theta}(\phi_1) \circ [\phi_1]) \left( [\phi_1]^{-1}(h) + \sum_{m \geq 1} (\tau_m^{\theta}(\phi_2) \circ [\phi_1]^{-1})(h) \right)}\\
           & = \displaystyle{-h + \kappa^{\theta}(\phi_1)\left( h + \sum_{m \geq 1} ([\phi_1]^{\otimes m} \circ \tau_m^{\theta}(\phi_2) \circ [\phi_1]^{-1})(h) \right) }\\
           & = -h + \kappa^{\theta}(\phi_1)(h) \\
           & \;\;\;\;\; \;\;\; + \kappa^{\theta}(\phi_1)( ([\phi_1]^{\otimes 2} \circ \tau_1^{\theta}(\phi_2) \circ [\phi_1]^{-1})(h))\\
           & \;\;\;\;\;  \;\;\; + \kappa^{\theta}(\phi_1)( ([\phi_1]^{\otimes 3} \circ \tau_2^{\theta}(\phi_2) \circ [\phi_1]^{-1})(h)) \;\; {\rm mod}\;  \widehat{U}_4.
 \end{array}
 $$
We note that
$$ \kappa^{\theta}(\phi)|_{H^{\otimes m}} = ({\rm id}_H + \tau^{\theta}(\phi))^{\otimes m} : H^{\otimes m} \longrightarrow H \times \widehat{U}_{2m}$$
for any $\phi \in {\rm Aut}(F)$ and so we have the following congruences mod $\widehat{U}_4$:
$$ \begin{array}{l}
\kappa^{\theta}(\phi_1)(h) \equiv h + \tau_1^{\theta}(\phi_1)(h) + \tau_2^{\theta}(\phi_1)(h),\\
\kappa^{\theta}(\phi_1)(([\phi_1]^{\otimes 2} \circ \tau_1^{\theta}(\phi_2) \circ [\phi_1]^{-1})(h)) \\
 \;\;\;\;\;\;\;\; \;\; \equiv  ([\phi_1]^{\otimes 2} \circ \tau_1^{\theta}(\phi_2) \circ [\phi_1]^{-1})(h)\\
 \;\;\;\;\;\;\;\;\;\;\;\;                     +  ( (\tau_1^{\theta}(\phi_1)\otimes {\rm id}_H + {\rm id}_H \otimes \tau_1^{\theta}(\phi_1)) \circ [\phi_1]^{\otimes 2} \circ \tau_1^{\theta}(\phi_2) \circ [\phi_1]^{-1})(h),\\
 \kappa^{\theta}(\phi_1)( ([\phi_1]^{\otimes 3} \circ \tau_1^{\theta}(\phi_2) \circ [\phi_1]^{-1})(h))  \equiv    ([\phi_1]^{\otimes 3} \circ \tau_2^{\theta}(\phi_2) \circ [\phi_1]^{-1})(h).
\end{array}
$$
Therefore we have
$$ \begin{array}{l} \tau^{\theta}(\phi_1\phi_2)(h)   \\
\equiv  \tau_1^{\theta}(\phi_1)(h)  +   \tau_2^{\theta}(\phi_1)(h)\\
\;\; +  ([\phi_1]^{\otimes 2} \circ \tau_1^{\theta}(\phi_2) \circ [\phi_1]^{-1})(h) \\
 \;\; +   ( (\tau_1^{\theta}(\phi_1)\otimes {\rm id}_H + {\rm id}_H \otimes \tau_1^{\theta}(\phi_1)) \circ [\phi_1]^{\otimes 2} \circ \tau_1^{\theta}(\phi_2) \circ [\phi_1]^{-1})(h)\\
\;\;  +     ([\phi_1]^{\otimes 3} \circ \tau_2^{\theta}(\phi_2) \circ [\phi_1]^{-1})(h) \;\;\;\;  \mbox{mod} \; \widehat{U}_4.
\end{array}
\leqno{(4.1.16.2)}
$$                                                                                           
Comparing (4.1.16.1) and (4.1.16.2), we obtain the assertions. $\;\;\Box$\\
\\
Next, we compute the $p$-Johnson maps for inner automorphisms of $F$.\\
\\
{\bf Proposition 4.1.17.}  {\it Let $ f \in F$ and $h \in H$.  For $m \geq 1$, we have}
$$ \tau_m^{\theta}({\rm Inn}(f))(h) = \theta_{m}(f)h + \sum_{j=1}^{m} \sum_{{\scriptstyle q_0+\cdots + q_j = m} \atop
{q_0 \geq 0, q_1,\dots,q_j \geq 1}} (-1)^j \theta_{q_0}(f)h\theta_{q_1}(f)\cdots \theta_{q_j}(f).$$                    
{\it In particular, we have, for $m=1, 2$,}
$$ \begin{array}{l}
\tau_1^{\theta}({\rm Inn}(f))(h) = [f] h - h [f],\\
\tau_2^{\theta}({\rm Inn}(f))(h) = \theta_2(f) h - h \theta_2(f) + h [f] [f] - [f] h [f].
\end{array}
$$
{\it Proof.}  Since $[{\rm Im}(f)] = {\rm id}_H$, by (4.1.12), we have 
$$ \begin{array}{ll}
\kappa^{\theta}({\rm Inn}(f))(z) & = (\theta \circ \widehat{{\rm Inn}(f)}\circ \theta^{-1})(z)\\
                                                       & = \theta(f) z \theta(f^{-1}) \\
                                                       & = \displaystyle{ \left(1+\sum_{m \geq 1} \theta_m(f)\right)z \left(1 + \sum_{j \geq 1}(-1)^j(\sum_{q \geq 1} \theta_q(f))^j \right)}
\end{array}
$$
for $z \in \widehat{U}$. Therefore, by (4.1.13), we have
$$ \begin{array}{ll}
\tau^{\theta}({\rm Inn}(f))(h) & = \kappa^{\theta}({\rm Im}(f))(h) - h \\
& = \displaystyle{ \sum_{m \geq 1} \theta_m(f)h + \sum_{j \geq 1} \sum_{{\scriptstyle q_o \geq 0.} \atop {\scriptstyle q_1,\dots,q_j \geq 1}} (-1)^j \theta_{q_0}(f)h\theta_{q_1}(f)\cdots \theta_{q_j}(f)}
\end{array}
$$
for $h \in H$. Taking the component in $H^{\otimes (m+1)}$, we obtain the assertion. $\;\; \Box$    \\                                              
\\                                                  
Finally we give the relation between the $p$-Johnson maps and  the $p$-Johnson homomorphisms in Section 2. \\
\\
{\bf Proposition 4.1.18.}  {\it The restriction of $\tau^{\theta}_m$ to  ${\rm A}_F(m)$  coincides with $\theta_{m+1}\circ \tau_m$ for each $m \geq 1$}:
$$ \tau^{\theta}_m|_{{\rm A}_F(m)} = \theta_{m+1} \circ \tau_m\, : \, {\rm A}_F(m) \longrightarrow {\rm Hom}(H, H^{\otimes (m+1)}),$$
where $\theta_{m+1}$ is the injection ${\rm gr}_{m+1}(F) \hookrightarrow H^{\otimes (m+1)}$ in (4.1.1).\\
\\
{\it Proof.} It suffices to show that for $\phi \in {\rm A}_F(m)$, 
$$ \tau_m^{\theta}(\phi)(X_j) = \theta_{m+1}(\tau_m(\phi)(X_j)) \;\; 1\leq j \leq r.$$
By (4.1.13) and $[\phi] = {\rm id}_H$, we have
$$ \begin{array}{ll}
\tau^{\theta}(\phi)(X_j) & = (\kappa^{\theta}(\phi)|_H - {\rm id}_H)(X_j)\\
 & = (\theta \circ \widehat{\phi} \circ \theta^{-1})(\theta(x_j)-1) - (\theta(x_j)-1)\\
 & = \theta(\phi(x_j)) - \theta(x_j).
 \end{array}
 $$
 Therefore we see that
 $$\tau_m^{\theta}(\phi)(X_j) = \mbox{ the component in} \; H^{\otimes (m+1)} \; \mbox{of}\; \theta(\phi(x_j)) - \theta(x_j).  \leqno{(4.1.18.1)}$$
 On the other hand, since $\phi(x_j)x_j^{-1} \in F_{m+1}$, we have
 $$ \theta(\phi(x_j)x_j^{-1}) \equiv 1 + \theta_{m+1}(\phi(x_j) x_j^{-1}) = 1 + \theta_{m+1}(\tau_m(\phi)(X_j)) \;\; {\rm mod}\; \widehat{U}_{m+2}.$$
 Multiplying the above equation by  $\theta(x_j)$ from right, we have
 $$ \theta(\phi(x_j)) \equiv \theta(x_j) + \theta_{m+1}(\tau_m(\phi)(X_j)) \;\; {\rm mod}\; \widehat{U}_{m+2}. \leqno{(4.1.18.2)}$$
By (4.1.18.1) and (4.1.18.2), we obtain the assertion. $\;\; \Box$ 
\\
\\
{\bf 4.2.  Examples in Non-Abelian Iwasawa theory.}  Let us come back to the arithmetic situation set up in 3.1 and keep the same notations.  So, as in (3.1.1) and (3.1.2),  we have an exact sequence of pro-$p$ Galois groups 
$$1 \longrightarrow G \longrightarrow {\cal G} \longrightarrow \Gamma \longrightarrow 1,$$
where 
$$G = {\rm Gal}(M/k_{\infty}), \; {\cal G} = {\rm Gal}(M/k) \; {\rm and}\; G = {\rm Gal}(k_{\infty}/k).
$$
In order to apply the materials in 4.1, we assume that 
$$
\mbox{(F)}\;\;\;\;\; G = {\rm Gal}(M/k_{\infty}) \; \mbox{is a  free pro-}p \; \mbox{group} \; F \; \mbox{on}\;  x_1,\dots , x_r.
$$
This condition (F) is satisfied for the following cases.
\\
\\
{\bf Example 4.2.1} ([Iw2], [W1]). Suppose that\\
(1) $k$ is totally real,\\
(2) $M := \tilde{k}_S$,\\
(3)  the Iwasawa $\mu$-invariant of $H_{\infty} = G/[G,G]$ is zero.\\
Then the condition (F) is satisfied where the generator rank $r$ is equal to the Iwasawa $\lambda$-invariant of $H_{\infty}$.\\
\\
To give the following example, we introduce the notation. For a field $K$, $K(p)$ denotes the maximal pro-$p$ extension of $K$.\\
\\
{\bf Example 4.2.2} ([Sc], [W2]). Suppose that \\
(1)  $k$ is a CM-field containing $\mu_p$ i.e., $k = k^+(\mu_p)$ where $k^+$ is the maximal totally real subfield of $k$,\\
(2)  the completions $k^+_{\frak{p}}$ of $k^+$ with respect to any prime $\frak{p}$ lying over $p$ do not contain $\mu_p$, \\
(3) $k_\infty$ is the cyclotomic $\mathbb{Z}_p$-extension of $k$, and\\
(4)  the Iwasawa $\mu$-invariant of the maximal Abelian unramified pro-$p$ Galois group over $k_\infty$ is zero.\\
The condition (4) is known to be true if $k$ is Abelian over $\mathbb{Q}$ ([FW]). So,  the above four conditions are satisfied for the $p$-th cyclotomic field $k= \mathbb{Q}(\mu_p)$, for instance.

A finite $p$-extension $L/k$ is called {\it positively ramified} over $S_p$ if $L_{\frak{p}} \subset k^+_{\frak{p}}(p)(\mu_p)$
for any prime $\frak{p}$ over $p$. Since the composite of positively ramified $p$-extensions is positively ramified again, the maximal 
positively ramified pro-$p$ extension of $k$ exists, and it contains the cyclotomic $\mathbb{Z}_p$-extension $k_{\infty}$. 
We then let

$M$ := the maximal pro-$p$ extension of $k$ which is unramified outside $S$ \\
$\;\;\;\;\;\;\;\;\;\; $ and positively ramified over $S_p$.\\
Then the condition (F) is satisfied with 
$$ r = 2\lambda^- + \#(S(k_{\infty}) \setminus S_p(k_{\infty})) - 1,$$
where $\lambda^-$ denotes the Iwasawa $\lambda^-$-invariant of $k$, and $S(k_{\infty})$ (resp. $S_p(k_{\infty})$) denotes the set of primes of $k_{\infty}$ lying over $S$ (resp. $S_p$). The pro-$p$ Galois group $F = G = {\rm Gal}(M/k_{\infty})$
has the following presentation
$$ \begin{array}{ll} F = & \langle a_1,b_1,\dots,a_{\lambda^-},b_{\lambda^-}, c_v (v \in S(k_{\infty}) \setminus S_p(k_{\infty})) \; | \\
                                      & \;\;\;\;\;\;\;\;\;\; \displaystyle{\prod_{v \in S(k_{\infty}) \setminus S_p(k_{\infty})} c_v \prod_{i=1}^{2\lambda^-}[a_i,b_i] = 1 \rangle}.
\end{array}$$
We may take $S$ to be $S_p \cup \{ \frak{q} \}$ such that there is only one prime of $k_{\infty}$ lying over $\frak{q}$ (there are infinitely many such $\frak{q}$). Then $F$ and ${\cal G}$ may be seen as arithmetic analogues of the fundamental groups of a one-boundary surface and a surface bundle over a circle (a fibered knot complement), respectively. \\

We fix a lift $\tilde{\gamma} \in {\cal  G}$ of a topological generator $\gamma$ of $\Gamma$ and consider the automorphism $\phi_{\tilde{\gamma}} := {\rm Inn}(\tilde{\gamma}) \in {\rm Aut}(F)$ as in (3.1.3). The $p$-power iterated action of $[\phi_{\tilde{\gamma}}]_m$ on $F/F_{m+1}$ is described by the $m$-th $p$-Johnson map
$$  \tau_m^{\theta} \, : \, {\rm Aut}(F) \longrightarrow {\rm Hom}_{\mathbb{F}_p}(H, H^{\otimes (m+1)}) \;\; (m \geq 1).$$
For an integer $d \geq 0$, we can  write 
$$ \tau_m^{\theta}(\phi_{\tilde{\gamma}}^{p^d})([f]) = \sum_{1 \leq i_1,\dots ,i_{m+1} \leq r} \tau^{\theta}(\phi_{\tilde{\gamma}}^{p^d})(i_1\cdots i_{m+1}; [f]) X_{i_1}\cdots X_{i_{m+1}}.$$
Suppose that $\phi_{\tilde{\gamma}}^{p^d} \in {\rm A}_F(m)$. Then we can also write
$$ \theta_{m+1} \circ \tau_m(\phi_{\tilde{\gamma}}^{p^d})([f]) = \sum_{1 \leq i_1,\dots ,i_{m+1} \leq r} \tau(\phi_{\tilde{\gamma}}^{p^d})(i_1\cdots i_{m+1}; [f]) X_{i_1}\cdots X_{i_{m+1}} \leqno{(4.2.3)} $$
and, by Proposition 4.1.18, we have 
$$\tau^{\theta}(\phi_{\tilde{\gamma}}^{p^d})(i_1\cdots i_{m+1}; X_j) = \tau(\phi_{\tilde{\gamma}}^{p^d})(i_1\cdots i_{m+1}; X_j) \in \mathbb{F}_p.  $$ 
These coefficients are numerical datum encoded in the Johnson maps/ homomorphisms. In Section 5, we express these coefficients in terms of Massey products in
Galois cohomology.\\

\begin{center}
{\bf 5. Massey products}
\end{center}

In this section, we give a cohomological interpretation of $p$-Johnson homomorphisms in terms of Massey products  in Galois cohomology.

A fixed prime number $p$ is arbitrary in 5.1 and assumed to be odd in 5.2.\\
\\
{\bf 5.1. Massey products and the Magnus expansion.} Firstly, we recall some general materials on Massey products. For the sign convention, we follow [Dw].   Let ${\cal G}$ be a pro-$p$ group and let  $\alpha_1,\dots,\alpha_m \in H^1({\cal G}, \mathbb{F}_p)$. A {\it Massey products}  $\langle \alpha_1,\dots,\alpha_m \rangle$ is said to be {\it defined} if there is an array
$$ A = \{a_{ij} \in C^1({\cal G},\mathbb{F}_p) \; | \; 1 \leq i \leq m+1, (i,j) \neq (1,m+1) \}$$
such that
$$\left\{ 
\begin{array}{l}
[a_{i,i+1}] = \alpha_i \;\; (1\leq i \leq m),\\
 \displaystyle{da_{ij} = \sum_{l=1}^{j-1} a_{1l}\cup a_{lj}} \;\; (j \neq i+1),

\end{array}
\right.
$$
where $d$ denotes the differential on cochains and $\cup$ denotes the cup product. An array $A$ is called a {\it defining system} for $\langle \alpha_1,\dots, \alpha_m \rangle$. Then we define $\langle \alpha_1,\dots, \alpha_m \rangle_A $ by the cohomology class represented by the $2$-cocycle 
$$  \sum_{l=2}^m a_{1l} \cup a_{l,m+1}.$$
A Massey product of $\alpha_1,\dots,\alpha_m$ is then defined by
$$ \langle \alpha_1,\dots, \alpha_m \rangle := \{ \langle \alpha_1,\dots, \alpha_m \rangle_A \in H^2({\cal G},\mathbb{F}_p) \; | \; A \; \mbox{ranges over defining systems} \}.$$

We recall some basic properties of Massey products, which will be used in 5.2.\\
\\
5.1.1.   One has $\langle \alpha_1, \alpha_2 \rangle = \alpha_1 \cup \alpha_2$. For $m \geq 3$, $ \langle \alpha_1,\dots, \alpha_m \rangle$ is defined and consists of a single element if 
$\langle \alpha_{i_1},\dots, \alpha_{i_l} \rangle = 0$ for all proper subsets $\{i_1,\dots, i_l \}$ of $\{ 1,\dots ,m\}$.\\
\\
5.1.2. Let $\Psi : {\cal G} \rightarrow {\cal G'}$ be a continuous  homomorphism of pro-$p$ groups. Then if $\langle \alpha_1,\dots,\alpha_m \rangle$ is defined for $\alpha_i \in H^1({\cal G'},\mathbb{F}_p)$
with defining system $A = (a_{ij})$, then so is $\langle \Psi^*(\alpha_1),\dots, \Psi^*(\alpha_m) \rangle$ with defining system $A^* = (\Psi^*(a_{ij}))$ and we have $\Psi^*(\langle \alpha_1,\dots,\alpha_m \rangle ) \subset \langle \Psi^*(\alpha_1),\dots ,\Psi^*(\alpha_m) \rangle$.\\

Next, we recall a relation between Massey products and the Magnus expansion. Let ${\cal G}$ be a finitely generated pro-$p$ group with a minimal presentation 
$$ 1 \longrightarrow N \longrightarrow {\cal F} \stackrel{\pi}{\longrightarrow} {\cal G} \longrightarrow 1,$$
where ${\cal F}$ is a free pro-$p$ group on $x_1, \dots, x_s$ with $s = {\rm dim}_{\mathbb{F}_p} H^1({\cal G}, \mathbb{F}_p)$. We set $g_i := \pi(x_i) \;(1\leq i \leq s)$. Note that $\pi$ induces the isomorphism 
$H^1({\cal G},\mathbb{F}_p) \simeq H^1({\cal F},\mathbb{F}_p).$ We let ${\rm tg} : H^1(N, \mathbb{F}_p)^{\cal G} \rightarrow H^2({\cal G},\mathbb{F}_p)$
 be the transgression map defined as follows. For $a \in H^1(N,\mathbb{F}_p)^{\cal G}$, choose a 1-cochain $b \in C^1({\cal F},\mathbb{F}_p)$ such that $b|_N = a$. Since the value $db(f_1,f_2)$, $f_i \in {\cal F}$, depends only on the cosets $f_i$ mod $N$, there is a 2-cocyle $c \in Z^2({\cal G},\mathbb{F}_p)$ such that $\pi^*(c) = db$. Then we define ${\rm tg}(a)$ by the class of $c$. By Hochschild-Serre spectral sequence,  ${\rm tg}$ is an isomorphism and so we have the dual isomorphism, called the Hopf isomorphism,
$${\rm tg}^{\vee} : H_2({\cal G},\mathbb{F}_p) \stackrel{\sim}{\rightarrow} H_1(N,\mathbb{F}_p)_{\cal G} = N/N^p[N,{\cal F}]. \leqno{(5.1.3)}$$
Then we have the following Proposition. The proof goes in the same manner as in [Ms1, Theorem 2.2.2].\\
\\
{\bf Proposition 5.1.4.}  {\it Notations being as above, let $\alpha_1,\dots,\alpha_m \in H^1({\cal G},\mathbb{F}_p)$ and  $A = (a_{ij})$ a defining system for the Massey product $\langle \alpha_1,\dots,\alpha_m \rangle$. Let $ f \in N$ and set $ \beta := ({\rm tg}^{\vee})^{-1}(f \; {\rm mod}\;  N^p[N,{\cal F}])$. Then we have} 
$$ \begin{array}{l} \langle \alpha_1,\dots, \alpha_m \rangle_A(\beta)\\
 = \displaystyle{\sum_{j=1}^m (-1)^{ j+1} \sum_{c_1+\cdots + c_j=m} \sum_{1\leq i_1,\dots,i_j \leq s}a_{1,1+c_1}(g_{i_1})\cdots a_{m+1-c_j,m+1}(g_{i_j})\epsilon(i_1\cdots i_j; f),}
\end{array}$$
{\it where  $c_1,\dots,c_j$ run over positive integers satisfying $c_1+\cdots + c_j = m$ and $g_i := \pi(x_i)$ $(1\leq i \leq s)$ and $\epsilon(i_1\cdots i_j; f)$ is the Magnus coefficient defined in $(4.1.9)$.}\\
\\
{\bf 5.2. Massey products and $p$-Johnson homomorphisms.} We come back to the arithmetic situation in 4.2 and keep the same notations. So we have an exact sequence of pro-$p$ Galois groups 
$$ 1 \longrightarrow F \longrightarrow {\cal G} \longrightarrow \Gamma \longrightarrow 1,$$
where 
$$F = {\rm Gal}(M/k_{\infty}), \; {\cal G} = {\rm Gal}(M/k) \; {\rm and}\; \Gamma = {\rm Gal}(k_{\infty}/k),$$
and $F$ is a free pro-$p$ group on $x_1,\dots,x_r$. We fix a lift $\tilde{\gamma} \in {\cal G}$ of a topological generator $\gamma$ of $\Gamma$ and let $\phi_{\tilde{\gamma}} := {\rm Inn}(\tilde{\gamma}) \in {\rm Aut}(F)$. 

Let $d(1)$ be the $p$-period of $[\phi_\gamma]$ on $H$ as in (3.2.3) so that $\phi_{\tilde{\gamma}}^{p^{d(1)}} \in {\rm IA}(F)$. If necessary, we replace the base field $k$ by the subextension $k_{d(1)}$ of $k_{\infty}$ with degree $[k_{d(1)}:k] =p^{d(1)}$ and $\gamma^{p^{d(1)}}$ with $\gamma$ so that we may suppose that
$$ \phi_{\tilde{\gamma}} \in {\rm IA}(F),$$
namely, $\phi_{\tilde{\gamma}}$ acts trivially on $H$.

For each integer $d \geq 0$, let $k_d$ be the subextension of $k_{\infty}$ with $[k_d : k] = p^d$ and let
$$ {\cal G}_d := {\rm Gal}(M/k_d).$$
Then the pro-$p$ group ${\cal G}_d$ has the presentation
$$ 1 \longrightarrow N_d \longrightarrow {\cal F} \stackrel{\pi_d}{\longrightarrow} {\cal G}_d \longrightarrow 1$$
where ${\cal F}$ is the free pro-$p$ group on $x_1,\dots , x_r, x_{r+1}$ with $\pi_d(x_{r+1}) = \gamma^{p^d}$ and $N_d$ is the closed subgroup of ${\cal F}$ generated normally by
$$ R_{j,d} :=  \phi_{\tilde{\gamma}}^{p^d}(x_j) (x_{r+1}x_j x_{r+1}^{-1})^{-1} \;\; (1\leq j \leq r).$$
\\
{\bf Lemma 5.2.1.} {\it For each integer $d\geq 0$, the homomorphism $\pi_d : {\cal F} \rightarrow {\cal G}_d$ induces the isomorphism of cohomology groups}
$$ H^1({\cal G}_d,\mathbb{F}_p) \stackrel{\sim}{\longrightarrow} H^1({\cal F},\mathbb{F}_p).$$
{\it Proof.} Since ${\cal G}_d = {\cal F}/N_d$, we have
$$ 
H^1({\cal G}_d,\mathbb{F}_p)  = {\rm Hom}_{\rm c}({\cal G}_d/{\cal G}_d^p[{\cal G}_d,{\cal G}_d], \mathbb{F}_p) \simeq {\rm Hom}_{\rm c}({\cal F}/N_d{\cal F}^p[{\cal F},{\cal F}], \mathbb{F}_p),
\leqno{(5.2.1.1)}
 $$
 where ${\rm Hom}_{\rm c}$ stands for the group of continuous homomorphisms. Since $\phi_{\tilde{\gamma}}^{p^d}$ acts trivially on $H = F/F^p[F,F]$, $\phi_{\tilde{\gamma}}^{p^d}(x_j) x_j^{-1}
 \in F^p[F,F]$ and so $R_{j,d} = \phi_{\tilde{\gamma}}^{p^d}(x_j)x_j^{-1} [x_j,x_{r+1}] \in {\cal F}^p[{\cal F},{\cal F}]$  $(1 \leq  j  \leq r)$. Therefore we have
 $$ N_d \subset {\cal F}^p[{\cal F},{\cal F}]. \leqno{(5.2.1.2)}$$
 By (5.2.1.1) and (5.2.1.2), we have
 $$ H^1({\cal G}_d,\mathbb{F}_p) \simeq {\rm Hom}_{\rm c}({\cal F}/{\cal F}^p[{\cal F},{\cal F}], \mathbb{F}_p) = H^1({\cal F},\mathbb{F}_p).  \;\; \Box$$
By Lemma 5.2.1, Hochschild-Serre spectral sequence yields the Hopf isomorphism as in (5.1.3)
$$ {\rm tg}^{\vee} \; : \; H_2({\cal G}_d,\mathbb{F}_p) \stackrel{\sim}{\longrightarrow} H_1(N_d,\mathbb{F}_p)_{{\cal G}_d} = N_d/N_d^p[N_d,{\cal F}],$$
and we define $\xi_{j,d} \in H_2({\cal G}_d,\mathbb{F}_p)$ by
$$ \xi_{j,d} := ({\rm tg}^{\vee})^{-1}(R_{j,d} \; {\rm mod} \; N_d^p[N_d,{\cal F}]) \;\; (1\leq j \leq r).$$

We set $g_j := \pi_d(x_i)$ $(1 \leq j \leq r+1)$ and let $g_i^* \in H^1({\cal G}_d,\mathbb{F}_p)$ denote the Kronecker dual to $g_j$, namely $g_i^*(g_j) = \delta_{ij}$.

For $d \geq 0$, let $m(d)$ be the integer defined in (3.2.5). Since $\phi_{\tilde{\gamma}} \in {\rm IA}(F)$, $m(d) \geq 1$. Let  $\tau_{m(d)}(\phi_{\tilde{\gamma}}^{p^d})(i_1\cdots i_{m(d)}; X_j)$ be the coefficients
of the $m(d)$-th  $p$-Johnson homomorphism defined  in (4.2.3). The following theorem gives an interpretation of $\tau_{m(d)}(\phi_{\tilde{\gamma}}^{p^d})(i_1\cdots i_{m(d)}; X_j)$ in terms of the Massey product in the cohomology of ${\cal G}_d$. \\
\\
{\bf Theorem 5.2.2.}  {\it Notations being as above, let $i_1, \dots , i_{m(d)+1} \in \{ 1, \dots , r \}$. Then the Massey product $\langle g_{i_1}^*,\cdots, g_{i_{m(d)+1}}^* \rangle$ is uniquely defined   and we have, for each $d\geq 0$, }
$$ \tau_{m(d)}(\phi_{\tilde{\gamma}}^{p^d})(i_1\cdots i_{m(d)+1};X_j) = (-1)^{m(d)+1}\langle g_{i_1}^*,\cdots, g_{i_{m(d)+1}}^* \rangle(\xi_{j,d}).$$
{\it Proof.} Let ${\cal G}_d'$ be the pro-$p$ group given by the presentation
$$1 \longrightarrow N_d' \longrightarrow F \stackrel{\pi_d'}\longrightarrow {\cal G}_d' \longrightarrow 1, $$
where $N'$ is the closed subgroup of $F$ generated normally by
$$ R'_{j,d} := \phi_{\tilde{\gamma}}^{p^d}(x_j) x_j^{-1} \;\; ( 1 \leq j \leq r). $$
We set $g'_j := \pi'(x_j)$ $(1\leq j \leq r)$ and let ${g'_i}^*$ be the Kronecker dual to $g'_j$. 
As in Lemma 5.2.1, $\pi_d'$ induces the isomorphism ${\rm tg} : H^1({\cal G}_d',\mathbb{F}_p) \stackrel{\sim}{\rightarrow} H^1(N_d',\mathbb{F}_p)_{{\cal G}'}$ and so we have the Hopf isomorphism 
 $ {\rm tg}^{\vee} : H_2({\cal G}_d',\mathbb{F}_p) \stackrel{\sim}{\rightarrow} H_1(N_d',\mathbb{F}_p)$. We define $\xi_{j,d}' \in H_2({\cal G}_d',\mathbb{F}_p)$ by $({\rm tg}^{\vee})^{-1}(R'_{j,d} \; {\rm mod}\; N_d'^p[N_d',F])$. 
Since $\phi_{\tilde{\gamma}}^{p^d} \in {\rm A}_F(m(d))$, we note $R'_{j,d} \in F_{m(d)+1} \; (1\leq j \leq r)$. 

Suppose $m(d) \geq 2$.  By Proposition 5.1.4, if $l \leq m(d)$, we have
$$ \langle  {g'_{i_1}}^*,\dots, {g'_{i_l}}^*  \rangle_{A'}(\xi'_{j,d}) = 0 $$
for any $i_1,\dots , i_l \in \{1,\dots r \}$, $1\leq j \leq r$, and any defining system $A'$, because we have $\epsilon(i_1\cdots i_l; R_{j,d}') = 0$. Since $R_{j,d}'$'s generate $H_1(N_d',\mathbb{F}_p)_{{\cal G}_d'}$, 
$ \langle  {g'_{i_1}}^*,\dots, {g'_{i_l}}^*  \rangle = 0$ for any $i_1,\dots , i_l \in \{1,\dots r \}$. Therefore, by 5.1.1, the Massey product $\langle  {g'_{i_1}}^*,\dots, {g'}^*_{i_{m(d)+1}} \rangle$ is uniquely defined and, by Proposition 5.1.4 again, we have
$$ \begin{array}{ll}
\langle  {g'_{i_1}}^*,\dots, {g'}^*_{i_{m(d)+1}}\rangle(\xi_{j,d}') & = (-1)^{m(d)+1} \epsilon(i_1\cdots i_{m(d)+1}; R'_{j,d}) \\
 & = (-1)^{m(d)+1} \tau_{m(d)}(\phi_{\tilde{\gamma}}^{p^d})(i_1\cdots i_{m(d)+1}; X_j).
 \end{array}
 \leqno{(5.2.2.1)}
$$

We define the homomorphism 
$$ \Psi \, : \, {\cal G}_d \longrightarrow {\cal G}_d'$$
by
$$ \Psi(g_j) := g'_j \; (1\leq j \leq r),\;\; \Psi(g_{r+1}) := 1.$$
so that we have
$$ \xi'_{j,d} = \Psi_*(\xi_{j,d})\;\; {g'_i}^* = \Psi^*(g_i^*) \;\; (1\leq i,j \leq r).$$
Then our assertion follows from (5.2.2.1) and the naturality 5.1.2 of Massey products as follows:
$$ \begin{array}{ll}
\langle g_{i_1}^*, \dots , g_{i_{m(d)+1}}^* \rangle (\xi_{j,d}) & = \langle \Psi^*({g'_i}^*),\dots , \Psi^*({g'}^*_{i_{m(d)+1}})\ \rangle (\Psi_*(\xi_{j,d}'))\\
 & = \Psi^*(\langle {g'}^*_{i_1}, \dots , {g'}^*_{i_{m(d)+1}} \rangle  ) (\Psi_*(\xi_{j,d}'))\\
  & = \langle  {g'}^*_{i_1},\dots, {g'}^*_{i_{m(d)+1}}  \rangle (\xi_{j,d}')\\
 & = (-1)^{m(d)+1} \tau_{m(d)}(\phi_{\tilde{\gamma}}^{p^d})(i_1\cdots i_{m(d)+1}; X_j). \;\;\;\; \Box
   \end{array}
 $$
\\
{\bf Remark 5.2.3.} (1) Theorem 5.2.2 may be regarded as an arithmetic analogue in non-Abelian Iwasawa theory of Kitano's result ([Ki, Theorem 4.1]).\\
(2) For Massey products in cohomology of a pro-$p$ group, we also refer to [G], [MT1] and [MT2].\\
\\
{\it Acknowledgement.} We would like to thank Yasushi Mizusawa, Manabu Ozaki, Takuya Sakasai and Takao Satoh
 for helpful communication. We would also like to thank the referee for useful comments.\\

\begin{flushleft}
{\bf References}\\
{[A]} S. Andreadakis, On the automorphisms of free groups and free nilpotent groups, Proc. London Math. Soc. {\bf 15}, (1965), 239-268.\\
{[BH]} J. P. Buhler, D. Harvey, Irregular primes to 163 million,  Math. Comp.  {\bf 80}  (2011),  no. 276, 2435-2444.\\
{[CFKSV]} J. Coates, T. Fukaya, K. Kato, R. Sujatha, O. Venjakob, The ${\rm GL}_2$  main conjecture for elliptic curves without complex multiplication, Publ. Math. Inst. Hautes Etudes Sci. No. {\bf 101}, (2005), 163-208.\\
{[Da]} M. Day, Nilpotence invariants of automorphism groups, Lecture note. Available at {\small http://www.math.caltech.edu/~2010-11/1term/ma191a/}\\
{[De]} C. Deninger, A note on arithmetic topology and dynamical systems, Algebraic number theory and algebraic geometry,  
Contemp. Math., {\bf 300}, Amer. Math. Soc., Providence, RI, 2002,  99-114.\\
{[DDMS]} J. D. Dixon, M. P. F. du Sautoy, A. Mann, D. Segal, Analytic pro-p  groups, Second edition. Cambridge Studies in Advanced Mathematics, {\bf 61}, Cambridge University Press, Cambridge, 1999. \\
{[Dw]} W. G.  Dwyer, Homology, Massey products and maps between groups, J. Pure Appl. Algebra  {\bf 6},  (1975), no. 2, 177-190.\\ 
{[FW]}  B. Ferrero, L. Washington, The Iwasawa invariant $\mu_p$ vanishes for abelian number fields, Ann. of Math.  {\bf 109},  (1979), no. 2, 377-395. \\
{[G]} J. G\"{a}rtner, Higher Massey products in the cohomology of mild pro-p -groups, J. Algebra  {\bf 422},  (2015), 788-820.\\
{[Ih]} Y. Ihara, On Galois representations arising from towers of coverings of $\mathbb{P}^1 \setminus \{0,1,\infty \}$,  Invent. Math.  {\bf 86},  (1986),  no. 3, 427-459.\\
{[IS]} K. Iwasawa, C. Sims, Computation of invariants in the theory of cyclotomic fields, J. Math. Soc. Japan, {\bf 18}, No.1, (1966), 86-96.\\
{[Iw1]} K. Iwasawa, On $\mathbb{Z}_l$-extensions of algebraic number fields, Ann. of Math. (2)  {\bf 98},  (1973), 246-326.\\
{[Iw2]} K. Iwasawa, Riemann-Hurwitz formula and $p$-adic Galois representations for number fields,
Tohoku Math. J. (2)  {\bf 33},  (1981), no. 2, 263-288. \\
{[Kt]} K. Kato, Iwasawa theory and generalizations, Proceedings of the International Congress of Mathematicians, Madrid, 2006, Vol. I,  Eur. Math. Soc., 335-357.\\
{[J]} D. Johnson, An abelian quotient of the mapping class group ${\cal T}_g$, Math. Ann. {\bf 249}, (1980), 225-242.\\
{[Kw]} N. Kawazumi,  Cohomological aspects of Magnus expansions, arXiv:0505497 [math.GT], 2006.\\
{[Ki]} T. Kitano, Johnson's homomorphisms of subgroups of the mapping class group, the Magnus expansion and Massey higher products of mapping tori, Topology and its application, {\bf 69}, 1996, 165-172.\\
{[Ko]} H. Koch, Galoissche Theorie der $p$-Erweiterungen, Deutscher Verlag der Wiss., Springer, 1970.\\
{[Ma]} B. Mazur,  Remarks on the Alexander polynomial, unpublished note. Available at http://www.math.harvard.edu/~mazur/older.html\\
{[MT1]} J. Min\'{a}\v{c}, N.D. T\^{a}n, Triple Massey products and Galois theory, J. Eur. Math. Soc. (JEMS)  {\bf 19},  (2017),  no. 1, 255-284.\\
{[MT2]} J. Min\'{a}\v{c}, N.D. T\^{a}n, The kernel unipotent conjecture and the vanishing of Massey products for odd rigid fields, (with an appendix written by I. Efrat, J. Min\'{a}\v{c} and N.D. T\^{a}n), Adv. Math.  {\bf 273},  (2015), 242-270.\\
{[Ms1]} M. Morishita, Milnor invariants and Massey products for prime numbers, Compos. Math.  {\bf 140},  (2004),  no. 1, 69-83. \\
{[Ms2]} M. Morishita, Knots and Primes -- An Introduction to Arithmetic Topology, Universitext, Springer, 2011.\\
{[Mt]} S. Morita, Abelian quotients of subgroups of the mapping class group of surfaces, Duke Math. J. {\bf 70}, (1993), 699-726.\\
{[O]} M. Ozaki, Non-abelian Iwasawa theory of $\mathbb{Z}_p$-extensions, J. Reine Angew. Math. {\bf 602}, (2007), 59-94.\\
{[Sa]} T. Satoh, A survey of the Johnson homomorphisms of the automorphism groups of free groups and related topics, arXiv:1204.0876[math.AT], 2012.\\
{[Sc]} A. Schmidt, Positively ramified extensions of algebraic number fields, J. Reine Angew. Math. {\bf 458}, (1995), 93-126.  \\
{[Sh]} R. Sharifi, Massey products and ideal class groups, J. Reine Angew. Math.  {\bf 603},  (2007), 1-33.\\
{[Wa]} L. Washington, Introduction to cyclotomic fields, Second edition. Graduate Texts in Math. {\bf 83}, Springer-Verlag, New York, 1997.\\
{[Wi1]} K. Wingberg, Duality theorems for $\Gamma$-extensions of algebraic number fields, Compositio Math. {\bf 55}, (1985), 333-381.\\
{[Wi2]} K. Wingberg,  Positiv-zerlegte $p$-Erweiterungen algebraischer Zahlk\"{o}rper, J. Reine Angew. Math.  {\bf 357},  (1985), 193-204. 
\end{flushleft}

\noindent M. Morishita:\\
Graduate School of Mathematics, Kyushu University, 744, Motooka, Nishi-ku, Fukuoka, 819-0395, Japan.\\
e-mail: morisita@math.kyushu-u.ac.jp \\
\\
Y. Terashima:\\
Department of Mathematical and Computing Sciences, Tokyo Institute of Technology, 2-12-1 Oh-okayama, Meguro-ku, Tokyo 152-8551, Japan.\\
e-mail: tera@is.titech.ac.jp

\end{document}